\documentclass[10pt,a4paper, abstracton,fleqn]{scrartcl}

 \pdfoutput=1

\usepackage[latin1]{inputenc}

\usepackage{amsthm,dcolumn,graphicx,multicol,fancyhdr}
\usepackage{latexsym,subfig,graphicx}
\usepackage{rotating,amsfonts, color,amssymb, amsmath, amscd, array, enumerate,afterpage}

\usepackage{natbib}

\usepackage{hyperref}

\usepackage{comment}

\newtheorem{theorem}{Theorem}[section]
\newtheorem{prop}[theorem]{Proposition}
\newtheorem{definition}{Definition}[section]
\newtheorem{lemma}[theorem]{Lemma}
\newtheorem{cor}[theorem]{Corollary}

\theoremstyle{definition}

\newtheorem{rem}{Remark}[section]

\newtheorem{assumption}{{Assumption $\mathcal{A}$.\negthickspace}}
\newtheorem{case}{{Case $\mathcal{C}$.\negthickspace}}

\renewenvironment{proof}[1][\negthickspace]{\textbf{Proof\thickspace#1.} }{\
\rule{0.5em}{0.5em}}

\makeatletter\@addtoreset{equation}{section}\makeatother
\makeatletter\@addtoreset{table}{section}\makeatother
\makeatletter\@addtoreset{figure}{section}\makeatother

\renewcommand{\cite}{\citep}\rm

\newcommand{\RR}{\mathbb{R}}

\newcommand{\cred}{\color{red}}

\newcommand{\var}{\operatorname{var}}
\newcommand{\E}{\operatorname{E}}
\newcommand{\cov}{\operatorname{cov}}
\newcommand{\ls}{\leqslant}
\newcommand{\gs}{\geqslant}
\newcommand{\eps}{\epsilon}

\newcommand{\dto}{\overset{\mathcal{L}}{\longrightarrow}}

\newcommand{\pto}{\overset{P}{\longrightarrow}}
\newcommand{\deq}{\overset{\mathcal{L}}{=}}

\newcommand{\po}{\sideset{_p}{}{\hspace{-0.9pt}\mathop\mathbf{o}}}
\newcommand{\qo}{\sideset{_q}{}{\hspace{-0.9pt}\mathop\mathbf{o}}}

\newcommand{\vth}{\vartheta}

\newcommand{\bp}{\mathbf{p}}
\newcommand{\bo}{\mathbf{o}}
\newcommand{\ba}{\mathbf{a}}
\newcommand{\be}{\mathbf{e}}
\newcommand{\br}{\mathbf{r}}
\newcommand{\bs}{\mathbf{s}}
\newcommand{\bM}{\mathbf{M}}

\newcommand{\bPhi}{\boldsymbol{\Phi}}
\newcommand{\bD}{\boldsymbol{\Delta}}

\newcommand{\bX}{\mathbf{X}}
\newcommand{\bZ}{\mathbf{Z}}
\newcommand{\HDE}{\mathcal{E}}

\newcommand{\hus}{Hu\v{s}kov\'{a} }

\newcommand{\hor}{Horv\'{a}th }

\newcommand{\haj}{H\'{a}jek }

\begin{document}

\title{ Efficiency of change point tests in high dimensional settings}

\author{ John A D Aston\footnote{Statistical Laboratory, DPMMS, University of Cambridge, Cambridge, CB3 9HD, UK; \texttt{j.aston@statslab.cam.ac.uk}} \and Claudia Kirch
\footnote{ Otto-von-Guericke University Magdeburg, Department of Mathematics, Institute of Mathematical Stochastics,
Postfach 4120, 39106 Magdeburg, Germany;
\texttt{claudia.kirch@ovgu.de}}}
\maketitle

\begin{abstract}
	\noindent
	While there is considerable work on change point analysis in univariate time series, more and more data being collected comes from high dimensional multivariate settings. This paper introduces the asymptotic concept of high dimensional efficiency which quantifies the detection power of different statistics in such situations. While being related to classic asymptotic relative efficiency, it is different in that it provides the rate at which the change can get smaller with dimension while still being detectable. This also allows for comparisons of different methods with different null asymptotics as is for example the case in high-dimensional change point settings.
Based on this new concept we investigate change point detection procedures using projections and develop asymptotic theory for how full panel (multivariate) tests compare with both oracle and random projections. Furthermore, for each given projection we can quantify a cone such that the corresponding projection statistic yields better power behavior if the true change direction is within this cone.
The effect of misspecification of the covariance on the power of the tests is investigated, because in many high dimensional situations estimation of the full dependency (covariance) between the multivariate observations in the panel is often either computationally or even theoretically infeasible. It turns out that the projection statistic is much more robust in this respect in terms of size and somewhat more robust in terms of power.  The theoretic quantification by the theory is accompanied by simulation results which confirm the theoretic (asymptotic) findings for surprisingly small samples. This shows in particular that the concept of high dimensional efficiency is indeed suitable to describe small sample power, and this is demonstrated in a multivariate example of market index data.

\end{abstract}

 {\bf Keywords: CUSUM; High Dimensional Efficiency; Model Misspecification; Panel Data; Projections}

 {\bf AMS Subject Classification 2000: 62M10; }

 \section{Introduction}

There has recently been a renaissance in research for statistical methods for change point problems \cite{HorvathR2014}. This has been driven by applications where non-stationarities in the data can often be best described as change points in the data generating process \cite{Eckley2010PFRK,frick2014multiscale,AK2}. However, data sets are now routinely considerably more complex than univariate time series classically studied in change point problems \cite{Page1954,RobbinsGLA2011,AueH2013,HorvathR2014}, and as such methodology for detecting and estimating change points in a wide variety of settings, such as multivariate \cite{horetal99,ombaoVSG2005,AueHHR2009,KirchMO2014} 
functional \cite{berkesetal09,aueetal09b,HormannK2010,AK1}
and high dimensional settings \cite{Bai2010,HorvathH2012,ChanHH2012,cho2015multiple} have recently been proposed. In panel data settings, these include methods based on taking maxima statistics across panels coordinate-wise \cite{jirak2015}, using sparsified binary segmentation for multiple change point detection \cite{cho2015multiple}, uses of double CUSUM procedures \cite{cho2015change}, as well as those based on structural assumptions such as sparsity \cite{WangS2016}.

Instead of looking at more and more complicated models, this paper uses a simple mean change setting to illustrate how the power is influenced in high dimensional settings. The results and techniques can subsequently be extended to more complex change point settings as well as different statistical frameworks, such as two sample tests. We make use of the following two key concepts: Firstly, we consider contiguous changes where the size of the change tends to zero as the sample size and with it the number of dimensions increases leading to the notion of high dimensional efficiency. This concept is closely related to Asymptotic Relative Efficiency (ARE) (see \citet[Sec. 3.4]{Lehmann1999} and \citet{lopesLW2011} where ARE is used in a high dimensional setting). Secondly, as a benchmark we investigate a class of tests based on projections, where the optimal (oracle) projection test is closely related to the likelihood ratio test under the knowledge of the direction of the change. Such tests can also be used in order to include a priori information about the expected change direction, where we can quantify how wrong the assumed direction can be and still yield better results than a full multivariate statistic which uses no information about the change direction. 

The aims of the paper are threefold: Firstly, we will investigate the asymptotic properties of tests based on projections as a plausible way to include prior information into the tests. Secondly, by using high dimensional efficiency, we consider several projection tests (including oracle and random projections as benchmarks) and compare them with the efficiency of existing tests that take the full covariance structure into account. Finally,  as in all high dimensional settings, the dependency between the components of the series can typically neither be effectively estimated nor even uniquely determined (for example if the sample size is less than the multivariate dimension) { unless restrictions on the covariance are enforced}. By considering the effect of misspecification of the model covariance on the size as well as efficiency, we can quantify the implications of this for different tests.

Somewhat obviously, highest efficiency can only be achieved under knowledge of the direction of the change. However, data practitioners, in many cases, explicitly have prior knowledge in which direction changes are likely to occur. It should be noted at this point, that changes in mean are equivalent to changes of direction in multivariate time series. In frequentist testing situations, practitioners' main interest is in test statistics which have power against a range of related alternatives while still controlling the size. For example, an economist may check the performance of several companies looking for changes caused by a recession. { There will often be a general idea as to which} sectors of the economy will gain or lose by the recession and therefore a good idea, at least qualitatively, as to what a change will approximately look like (downward resp. upward shift depending on which sector a particular company is in) if there is a change present. Similarly, in medical studies, it will often be known a-priori whether genes are likely to be co-regulated causing changes to be in similar directions for groups of genes in genetic time series.

Incorporating this a-priori information about how the change affects the components by using corresponding projections can lead to a considerable power improvement if the change is indeed in the expected direction. It is also important that, as in many cases the a-priori knowledge is qualitative, the test has higher power than standard tests not only for that particular direction but also for other directions close by. Additionally, these projections lead to tests where the size is better controlled if no change is present. Using the concept of high dimensional efficiency allows the specification of a cone around a given projection such that the projection statistic has better power than the multivariate/panel statistic if the true change is within this cone. In addition, while the prior information itself might be reliable, inherent misspecification in other parts of the model, such as the covariance structure, will have a detrimental effect on detection, and it is of interest to quantify the effect of these as well.

The results in this paper will be benchmarked against taking the simple approach of using a random projection in {a single direction} to reduce the dimension of the data. Random projections are becoming increasingly popular in high dimensional statistics with applications in Linear Discriminant Analysis \cite{durrantk2010} and two sample testing \cite{lopesLW2011,SrivastavaLR2014}. This is primarily based on the insight from the Johnson-Lindenstrauss lemma that an optimal projection in the sense that the distances are preserved for a given set of data is independent of the dimension of the data \cite{johnsonl1984} and thus random projections can often be a useful way to perform a dimension reduction for the data \cite{BaraniukDDW2008}.
However, in our context, we will see that a random projection will not work as well as truly multivariate methods, let alone projections with prior knowledge, but can only serve as a lower benchmark.

We will consider a simple setup for our analysis, although one which is inherently the base for most other procedures, and one which can easily be extended to complex time dependencies and change point definitions using corresponding results from the literature \cite{kirchTK2014report,kirchTK2014}. For a set of observations $X_{i,t}$, $ 1\ls i\ls d=d_T, 1\ls t\ls T$, the change point model is defined to be \begin{align}\label{eq_model}
	X_{i,t}=\mu_i+\delta_{i,T} \,g(t/T)+e_{i,t}, \quad 1\ls i\ls d=d_T, 1\ls t\ls T,
\end{align}
where $\E e_{i,t}=0$ for all $i$ and $t$ with $0<\sigma_i^2=\var e_{i,t}<\infty$ and $g:[0,1]\to \RR$ is a Riemann-integrable function. Here $\delta_{i,T}$ indicates the size of the change for each component. This setup incorporates a wide variety of possible changes by the suitable selection of the function $g$, as will be seen below. For simplicity, for now it is assumed that $\{e_{i,t}:t\in\mathbb{Z}\}$ are independent,  i.e.\ we assume independence across time but not location. If the number of dimensions $d$ is fixed, the results readily generalise to situations where a multivariate functional limit theorem exists as is the case for many weak dependent time series. If $d$ can increase to infinity with $T$, then generalizations are possible if the $\{e_{i,t}:1\ls t\ls T\}$ form a linear process in time but the errors are independent between components (dependency between components will be discussed in detail in the next section). Existence of moments strictly larger than two is needed in all cases. Furthermore, the developed theory applies equally to one- and two-sample testing and can be seen as somewhat analogous to methods for multivariate adaptive design \cite{MinasAS2014}. 

The change (direction) is given by $\bD_d=(\delta_{1,T},\ldots,\delta_{d,T})^T$ and the type of alternative is given by the function $g$ in rescaled time. While $g$ is defined in a general way, it includes as special cases most of the usual change point alternatives, for example,
\begin{itemize}
\item At most one change (AMOC): $g(u) = \left\{\begin{array}{ll}0& 0\leq u \leq \theta\\1& \theta < u \leq 1\end{array}\right.$
\item Epidemic change (AMOC): $g(u) = \left\{\begin{array}{ll}0& 0\leq u \leq \theta_1\\1& \theta_1 < u < \theta_2\\0& \theta_2 < u \leq 1\end{array}\right.$
\end{itemize}
The form of $g$ will influence the choice of test statistic to detect the change point. As in the above two examples in the typical definition of change points the function $g$ is modelled by a step function (which can approximate many smooth functions well). In such situations,  test statistics based on partial sums of the observations have been well studied \cite{csoehor}. It will be shown that statistics based on partial sums are robust (in the sense of still having non-zero power) to a wide variety of $g$.


The model in \eqref{eq_model} is defined for univariate ($d=1$), multivariate ($d$ fixed) or panel data ($d\to\infty$).  The panel data (also known as ``small n large p'' or ``high dimensional low sample size'') setting is able to capture the small sample properties very well in situations where $d$ is comparable or even larger than $T$ using asymptotic considerations. In this asymptotic framework the detection ability or efficiency of various tests can be defined by the rates at which vanishing alternatives can still be detected. However, many of our results, particularly for the proposed projection tests, are also qualitatively valid in the multivariate or $d$ fixed setting.





The paper proceeds as follows. In Section \ref{s:CPPs}, the concept of high dimensional efficiency as a way of comparing the power of high dimensional tests is introduced. This is done using projection statistics, which will also act as benchmarks. In Section \ref{section_comparison}, the projection based statistics will be compared with the panel based change point statistics already suggested in \citet{HorvathH2012}, both in terms of control of size and efficiency, particular with relation to the (mis)specification of the dependence structure. Section \ref{sec_data} provides a short illustrative example with respect to multivariate market index data. Section \ref{section_conclusions} concludes with some discussion of the different statistics proposed, while Section \ref{section_proofs} gives the proofs of the results in the paper. In addition, rather than a separate simulation section, simulations will be interspersed throughout the theory. They complement the theoretic results, confirming that the conclusion are already valid for small samples, thus verifying that the concept of high-dimensional efficiency is indeed suitable to understand the power behavior of different test statistics. In all cases the simulations are based on 1000 repetitions of i.i.d. normally distributed data for each set of situations, and unless otherwise stated the number of time points is $T=100$ with the change (if present) occurring half way through the series. Except in the simulations concerning size itself, all results are empirically size corrected to account for the size issues for the multivariate (panel) statistic that will be seen in Figure \ref{f:size}.

\section{Change Points and Projections}\label{s:CPPs}
\subsection{High dimensional efficiency}
As the main focus of this paper is to compare several test statistics with respect to their detection power, we introduce a new asymptotic concept that allows us to understand this detection power in a high dimensional context. In the subsequent sections, simulations accompanying the theoretic results will show that this concept is indeed able to give insight into the small sample detection power.

Consider a typical testing situation, where (possibly after reparametrization) we test
\begin{align}
	H_0: \mathbf{a}=0,\quad \text{ against }\quad H_1:\mathbf{a}\neq 0.\label{eq_test_HDE}
\end{align}
Typically, large enough alternatives, will be detected by all reasonable statistics for a given problem. In asymptotic theory this corresponds to fixed alternatives, where $a=c\neq 0$, for which tests typically have asymptotic power one.

To understand the small sample power of different statistics such asymptotics are therefore not suitable. Instead the asymptotics for local or contiguous alternatives with $a=a_T\to 0$ are considered.
For a panel setting we define:
\begin{definition}\label{def_HDE}
	Consider the testing situation \eqref{eq_test_HDE} with sample size $T\to\infty$ and sample dimension $d=d_T\to \infty$.  A test with statistic $\mathcal{T}(\mathbf{X}_1,\ldots,\mathbf{X}_T)$ has (absolute) {\bf high dimensional efficiency} $\HDE(\mathbf{a}_d) $ for a sequence of alternatives $\mathbf{a}_d$  if
	\begin{enumerate}[(i)]
		\item $\mathcal{T}(\mathbf{X}_1,\ldots,\mathbf{X}_T)\dto L$ for some non-degenerate limit distribution $L$ under $H_0$,
		\item $\mathcal{T}(\mathbf{X}_1,\ldots,\mathbf{X}_T)\pto \infty$ if $\sqrt{T}\,\HDE(\mathbf{a}_d)\to\infty$,
		\item $\mathcal{T}(\mathbf{X}_1,\ldots,\mathbf{X}_T)\dto L$ if $\sqrt{T}\,\HDE(\mathbf{a}_d)\to 0$.
	\end{enumerate}
Note that the $\HDE(\mathbf{a}_d)$ is only defined up to multiplicative constants, and has to be understood as a representative of a class.
\end{definition}
In particular this shows that the asymptotic power is one if $\sqrt{T}\,\HDE(\mathbf{a}_d)\to\infty$, but equal to the level if $\sqrt{T}\,\HDE(\mathbf{a}_d)\to 0$.
Typically, for $\sqrt{T}\,\HDE(\mathbf{a}_d)\to \alpha\neq 0$ it holds $\mathcal{T}(\mathbf{X}_1,\ldots,\mathbf{X}_T)\dto L(\alpha)\overset{\mathcal{D}}{\not =} L$, usually resulting in an asymptotic power strictly between the level and one.
In the  classic notion (with $d$ constant) of absolute relative efficiency (ARE, or Pitman Efficiency) for test statistics with a standard normal limit
it is the additive shift between $L(\alpha)$ and $L$ \cite[see][Sec 3.4]{Lehmann1999}) that shows power differences for different statistics. Consequently, this shift has been used to define asymptotic efficiency. For different null asymptotics the comparison becomes much more cumbersome as the quantiles of the different limit distributions had to be taken into account as well. In our definition above, we concentrate on the efficiency in terms of the asymptotic rates with respect to the increasing dimension (as two estimators are equivalent up to a dimension free constant). Should the rates be equivalent, classic notions of ARE then apply, although with the usual difficulties should the limit distributions be different.

In the asymptotic panel setup, on the other hand, the differences with respect to the dimension $d$ are now visible in the rates, with which contiguous alternatives can disappear and still be asymptotically detectable. Therefore, we chose this rate to define asymptotic high dimensional efficiency.
Additionally, it is no longer required that different test statistics have the same limit distribution under the null hypothesis (which would be a problem  in this paper).

\subsection{Projections}
We now describe how projections can be used to obtain change point statistics in high dimensional settings, which will be used as both benchmark statistics for a truly multivariate statistic as well as a reasonable alternative if some knowledge about the direction of the change is present.

In model \eqref{eq_model}, the change $\boldsymbol{\Delta}_d=(\delta_{1,T},\ldots,\delta_{d,T})^T$ (as a direction) is always a rank one (vector) object no matter the number of components $d$. This observation suggests that knowing the direction of the change $\boldsymbol{\Delta}_d$ in addition to the underlying covariance structure can significantly increase the signal-to-noise ratio. Furthermore, for $\mu$ and $\boldsymbol{\Delta}_d/\|\boldsymbol{\Delta}_d\|$ (but not $\|\boldsymbol{\Delta}_d\|$) known with i.i.d.\ normal errors, one can easily verify that the corresponding likelihood ratio statistic is obtained as a projection statistic with projection vector $\Sigma^{-1}\boldsymbol{\Delta}_d$, which can also be viewed as an oracle projection.
Under \eqref{eq_model} it holds
\begin{align*}
	\langle \mathbf{X}_d(t),\bp_d\rangle=\langle \boldsymbol{\mu},\bp_d\rangle+\langle \bD_d,\bp_d\rangle g(t/T) +\langle \be_t,\bp_d\rangle,
\end{align*}
where $\mathbf{X}_d(t)=(X_{1,t},\ldots,X_{d,T})^{T}$, $\boldsymbol{\mu}=(\mu_1,\ldots,\mu_d)^T$ and $\be_t=(e_{1,t},\ldots,e_{d,t})^T$. The projection vector $\bp_d$ plays a crucial role in the following analysis and will be called the search direction. This representation shows that the projected time series exhibits the same behavior as before as long as the change is not orthogonal to the projection vector. Furthermore, the power is the better the larger $\langle \bD_d,\bp_d\rangle$ and the smaller the variance of $\langle \be_t,\bp_d\rangle$ is. Consequently, an optimal projection in terms of power depends on $\bD_d$ as well as $\Sigma=\var \be_1$.
In applications, certain changes are either expected or of particular interest e.g.\ an economist looking at the performance of several companies expecting changes caused by a recession will have a good idea which companies will profit or lose. This knowledge can  be used to increase the power in directions close to the search direction $\bp_d$ while decreasing it for changes that are close to orthogonal to it. Using projections can furthermore robustify the size of the test under the null hypothesis with respect to misspecification and estimation error.

In order to qualify this informal statement, we will consider contiguous changes for several change point tests, where $\|\boldsymbol\Delta_d\|\to 0$ but with such a rate that the power of the corresponding test is strictly between the size and one, as indicated in the previous subsection. 

In order to be able to prove asymptotic results for change point statistics based on  projections even if $d\to\infty$, we need to make the following assumptions on the underlying error structure. This is much weaker than the independence assumption as considered by \citet{HorvathH2012}. Furthermore, we do not need to restrict the rate with which $d$ grows.
If we do have restrictions on the growth rate in particular for the multivariate setting with $d$ fixed, these assumptions can be relaxed and more general error sequences can be allowed.

\begin{assumption}\label{model_VAR0}Let $\eta_{1,t}(d),\eta_{2,t}(d),\ldots$ be independent with $\E \eta_{i,t}(d)=0$, $\var \eta_{i,t}(d)=1$ and $\E |\eta_{i,t}(d)|^{\nu}\ls C<\infty$ for some $\nu>2$ and all $i$ and $d$. For $t=1,\ldots,T$ we additionally assume for simplicity that $(\eta_{1,t}(d),\eta_{2,t}(d),\ldots)$ are identically distributed (leading to data which is identically distributed across time).
	The errors within the components are then given as linear processes of these innovations:
\begin{align*}
	e_{l,t}(d)=\sum_{j\gs 1}a_{l,j}(d)\eta_{j,t}(d), \quad l=1,\ldots,d,\quad \sum_{j\gs 1}a_{l,j}(d)^2<\infty
\end{align*}
or equivalently in vector notation $e_t(d)=(e_{1,t}(d),\ldots,e_{d,t}(d))^T$
 and $\boldsymbol a_j(d)=(a_{1,j}(d),\ldots,a_{d,j}(d))^T$
 \begin{align*}
	 \be_t(d)=\sum_{j\gs 1}\ba_j(d)\eta_{j,t}(d).
 \end{align*}
\end{assumption}

These assumptions allow us to considered many varied dependency relationships between the components (and we will concentrate on within the component dependency at this point, as temporal dependency adds multiple layers of notational difficulties, but little in the way of insight as almost all results generalise simply for weakly dependent and linear processes).

The following three  cases of different dependency structures are very helpful in understanding different effects that can occur and will be used as examples throughout the paper:
\begin{case}[Independent Components]\label{case_ind}
The components are independent, i.e.\
$\ba_j=(0,\ldots,s_j,\ldots,0)^T$ the vector which is $s_j> 0$ at point $j$ and zero everywhere else, $j\ls d$, and $\ba_j=\boldsymbol 0$ for $j\gs d+1$.  In particular, each channel has variance
			\begin{align*}
                \sigma_j^2=s_j^2.
            \end{align*}
\end{case}

\begin{case}[Fully Dependent Components]\label{case_dep}
	There is one common factor to all components, leading to  completely dependent components, i.e.\
$\ba_{1}=\bPhi_d=(\Phi_1,\ldots,\Phi_d)^T$, $\ba_j=\boldsymbol 0$ for $j\gs 2$. In this case,
			\begin{align*}
				\sigma_j^2=\Phi_j^2.
			\end{align*}
This case, while being somewhat pathological, is useful for gaining intuition into the effects of possible dependence and also helps with understanding the next case.		\end{case}

\begin{case}[Mixed Components]\label{case_mixed}
The components contain both an independent and dependent term. Let $\ba_j=(0,\ldots,s_j,\ldots,0)^T$ the vector which is $s_j>0$ at point $j$ and zero everywhere else, and $\ba_{d+1}=\bPhi_d=(\Phi_1,\ldots,\Phi_d)^T$, $\ba_j=\boldsymbol 0$ for $j\gs d+2$. Then
			\begin{align*}
				\sigma_j^2=s_j^2+\Phi_j^2
            \end{align*}
	    This mixed case allows consideration of dependency structures between cases $\mathcal{C.}$\ref{case_ind} and $\mathcal{C.}$\ref{case_dep}. It is used in the simulations with $\Phi_d=\Phi (1,\ldots,1)^T$, where $\Phi=0$ corresponds to $\mathcal{C}.1$ and $\Phi\to\infty$ corresponds to $\mathcal{C}.3$. We also use this particular example for the panel statistic in Section~\ref{section_panel_dep} to quantify the effect of misspecification.
\end{case}

Of course, many other dependency structures are possible, but these three cases give insight into the cases of no, complete and some dependency respectively. { In particular, as the change is always rank one, taking a rank one form of dependency, as in cases $\mathcal{C.}$\ref{case_dep} and as part of $\mathcal{C.}$\ref{case_mixed}, still allows somewhat general conclusions to be drawn.}

\subsection{Change point statistics}\label{sec_cps}
Standard statistics such as the CUSUM statistic are based on partial sum processes, so in order to quantify the possible power gain by the use of projections we will consider the partial sum process of the projections, i.e.\
\begin{align}
	&U_{d,T}(x)=\langle \bZ_T(x),\bp_d\rangle=\frac{1}{\sqrt{T}}\sum_{t=1}^{\lfloor Tx\rfloor}\left( \langle \bX_d(t), \bp_d\rangle -\frac{1}{T}\sum_{j=1}^T\langle \bX_d(j),\bp_d\rangle \right),\\
	&Z_{T,i}(x)=\frac{1}{T^{1/2}}\left( \sum_{t=1}^{\lfloor Tx\rfloor}X_{i,t}-\frac{\lfloor Tx\rfloor}{T}\sum_{t=1}^TX_{i,t} \right),\label{eq_def_Z}
\end{align}
where $\bX_d(t)=(X_{1,1},\ldots,X_{d,T})^T$.

Different test statistics can be defined for a range of $g$ in \eqref{eq_model}, however, assuming that $g \not \equiv 0$, the hypothesis of interest is
\[
H_0: \bD_{d}=\boldsymbol 0
\]
versus the alternative
\begin{align*}
H_1: \bD_{d} \neq \boldsymbol 0.
\end{align*}

Test statistics are now defined in order to give good power characteristics for a particular $g$ function. For example, the classic AMOC statistic for univariate and multivariate change point detection is based on $U_{d,T}(x)/\tau(\bp_d)$, with
\begin{align}\label{eq_tau}
	\tau^2(\bp_d)=\bp_d^T\var\left(\be_1(d) \right)\bp_d.
\end{align}
Typically, either the following max or sum type statistics are used:

\begin{align*}	\max_{1\ls k\ls T}w(k/T)\left|\frac{U_{d,T}(k/T)}{\tau(\bp_d)}\right|,\qquad \frac{1}{T}\sum_{k=1}^Tw(k/T)\left|\frac{U_{d,T}(k/T)}{\tau(\bp_d)}\right|,
\end{align*}
where $w\gs 0$ is continuous (which can be relaxed) and fulfills \eqref{eq_weights} (confer e.g.\ the book by \citet{csoehor}). The choice of weight function $w(\cdot)$ can increase power for certain locations of the change points \cite{KirchMO2014}.

 For the epidemic change, typical test statistics are given by
\begin{align*}
	&\max_{1\ls k_1<k_2\ls T}\frac{1}{\tau(\bp_d)}|U_{d,T}\left( k_2/T \right)-U_{d,T}\left( k_1/T \right)|,\qquad\\
	&  \text{or}\\
	&\frac{1}{T^2}\sum_{1\ls k_1<k_2\ls T}\frac{1}{\tau(\bp_d)}|U_{d,T}\left( k_2/T \right)-U_{d,T}\left( k_1/T \right)|.
\end{align*}

In the next section we first derive a functional central limit theorem for the process $U_{d,T}(x)$, which implies the asymptotic null behavior for the above tests. Then, we derive the asymptotic behavior of the partial sum process under contiguous alternatives  to obtain the high dimensional efficiency for projection statistics. 

\subsection{Efficiency of Change point tests based on projections}
In this section, we derive the efficiency of change point tests based on projections under rather general assumptions. Furthermore, we will see that the size behavior is very robust with respect to deviations from the assumed  underlying covariance structure. The power on the other hand turns out to be less robust but more so than statistics taking the full multivariate information into account.

\subsubsection{Null Asymptotics}
As a first step towards the efficiency of projection statistics, we derive the null asymptotics. This is also of independent interest if projection statistics are applied to a given data set in order to find appropriate critical values. In the following theorem $d$ can be fixed but it is also allowed that $d=d_T\to\infty$, where no restrictions on the rate of convergence are necessary.
\begin{theorem}\label{th_proj_null}
	Let model \eqref{eq_model} hold. Let $\bp_d$ be a possibly random projection independent of $\{e_{i,t}: 1\ls t\ls T, 1\ls i\ls d\}$.
 Furthermore, let $\bp_d^T\cov(\be_1(d)) \bp_d\neq 0$ (almost surely), which means that the projected data is not degenerate with probability one.
	\begin{enumerate}[a)]
		\item Under Assumption $\mathcal{A}$.\ref{model_VAR0} and if $\{\bp_d\}$ is independent of $\{\eta_{i,t}(d):i\gs 1,1\ls t\ls T\}$, then it holds under the null hypothesis
\begin{align}
	\left\{\frac{U_{d,T}(x)}{\tau(\bp_d)}:0\ls x\ls 1\,|\,\bp_d\right\}\overset{D[0,1]}{\longrightarrow} \{B(x):0\ls x\ls 1\}\qquad a.s.,\label{eq_indep_struc}
\end{align}
where $B(\cdot)$ is a standard Brownian bridge.
\item For i.i.d.\ error sequences $\{\be_t(d):t=1 ,\ldots,d\}$, $\be_t(d)=(e_{1,t}(d),\ldots,e_{d,t}(d))^T$ with an arbitrary dependency structure across components, and if $\E |e_{1,t}(d)|^{\nu}\ls C<\infty$ for all $t$ and $d$ as well as	
	\begin{align}\label{eq_general_dep}
		\frac{\|\bp_d\|_1^2}{\bp_d^{T}\cov (\be_t)\bp_d^T}=o(T^{1-2/\nu}) \quad a.s.,
	\end{align}
	where $\|\mathbf{a}\|_1=\sum_{j=1}^d|a_j|$,
	then \eqref{eq_indep_struc} holds.
	\end{enumerate}
	The assertions remain true if $\tau^2(\bp_d)$ is replaced by $\widehat{\tau}^2_{d,T}$ such that for all $\eps >0$
	\begin{align}\label{condvar}
		P\left(\left|\frac{\widehat{\tau}^2_{d,T}}{\tau^2(\bp_d)}- 1\right|>\eps\right) \to 0\qquad a.s.
	\end{align}
\end{theorem}
Assumption
\eqref{eq_general_dep} is always fulfilled for the multivariate situation with $d$ fixed or if $d$ is growing sufficiently slowly with respect to $T$ as the left hand side of \eqref{eq_general_dep} is always bounded by $\sqrt{d}$ if $\bp_d^T\cov(e)\bp_d/\|p_d\|^2$ is bounded away from zero. Otherwise, the assumption may hold for certain projections but not others. However, in this case, it is possible to put stronger assumptions on the error sequence such as in a), which are still much weaker than the usual assumption for panel data, that components are independent.
In these cases projection methods hold the size asymptotically, no matter what the dependency structure between components is and without having to estimate this dependency structure.

This is in contrast to the multivariate statistic which suffers from considerable size distortions if this underlying covariance structure is estimated incorrectly. The estimation of the covariance structure is a difficult problem in higher dimensions in particular since an estimator for the inverse is needed with additional numerical problems arising. The problem becomes even harder if time series errors are present, in which case the long-run covariance rather than the covariance matrix needs to be estimated \cite{HormannK2010,AK2,KirchMO2014}. While the size of the projection procedure is unaffected by the underlying dependency across components, we will see in the next section that for optimal efficiency hence power we need not only to know the change $\bD_d$ but also the inverse of the covariance matrix. Nevertheless the power of projection procedures turns out to be more robust with respect to misspecification than a size-corrected panel statistic, that takes the full multivariate information into account.

The following lemma shows that the following two different estimators for $\tau(\bp_d)$ under the null hypothesis are both consistent. The second one is typically still  consistent in the presence of one mean change which usually leads to a power improvement in the test for small samples. An analogous version can be defined  for the epidemic change situation. However, it is much harder to get an equivalent correction in the multivariate setting because the covariance matrix determines how different components are weighted, which in turn has an effect on the location of the maximum. This problem does not arise in the univariate situation, because the location of the maximum does not depend on the variance estimate.

\begin{lemma}\label{lem_variance}
Consider
\begin{align}\label{eq_est_tau}
	\widehat{\tau}^2_{1,d,T}(\bp_d) =\frac 1 T\sum_{j=1}^T\left(\bp_d^T\be_t(d)-\frac 1 T \sum_{i=1}^T\bp_d^T\be_t(d)\right)^2
\end{align}
as well as
\begin{align}\label{eq_est_tau2}
	&	\widehat{\tau}^2_{2,d,T}(\bp_d) =\frac 1 T\left(\sum_{j=1}^{\widehat{k}_{d,T}}\left(\bp_d^T\be_j(d)-\frac 1 T \sum_{i=1}^{\widehat{k}_{d,T}}\bp_d^T\be_i(d)\right)^2+\sum_{j=\widehat{k}_{d,T}+1}^{T}\left(\bp_d^T\be_t(d)-\frac 1 T \sum_{i=\widehat{k}_{d,T}+1}^T\bp_d^T\be_i(d)\right)^2\right), \\
	&\text{where}\quad \widehat{k}_{d,T}=\arg\max_{t=1,\ldots,T}U_{d,T}(t/T).\notag
\end{align}
\begin{enumerate}[a)]
	\item Under the assumptions of Theorem~\ref{th_proj_null} a) both estimators \eqref{eq_est_tau} as well as \eqref{eq_est_tau2} fulfill \eqref{condvar}.
	\item Under the assumptions of Theorem~\ref{th_proj_null} b), then estimator \eqref{eq_est_tau} fulfills \eqref{condvar} under the  assumption
		\begin{align*}		\frac{\|\bp_d\|_1^2}{\bp_d^{T}\cov (\be_t)\bp_d^T}=o(T^{1-2/\min(\nu,4)})\qquad a.s.,
	\end{align*}
	while estimator \eqref{eq_est_tau2} fulfills it under the assumption
	\begin{align*}		\frac{\|\bp_d\|_1^2}{\bp_d^{T}\cov (\be_t)\bp_d^T}=o(T^{1-2/\min(\nu,4)} (\log T)^{-1})\qquad a.s.,
	\end{align*}
\end{enumerate}
\end{lemma}

The following theorem gives the null asymptotic for the simple CUSUM statistic for the at most one change,  other statistics as given in Section~\ref{sec_cps} can be dealt with along the same lines.

\begin{cor}\label{th_asym_null_stat}
	Let the assumptions of Theorem~\ref{th_proj_null} be fulfilled and $\widehat{\tau}(\bp_d)$ fulfill \eqref{condvar}  under the null hypothesis, then for all $x\in\RR$ it holds under the null hypothesis
	\begin{align*}
		&P\left(	\max_{1\ls k\ls T}w^2(k/T) \frac{U_{d,T}^2(k/T)}{\widehat{\tau}^2(\bp_d)}\ls x\,\Big|\,\bp_d\right)\to P\left( \max_{0\ls t\ls 1}w^2(t)B^2(t) \ls x \right)\quad a.s.\\
		&P\left(\frac 1 T	\sum_{1\ls k\ls T}w^2(k/T) \frac{U_{d,T}^2(k/T)}{\widehat{\tau}^2(\bp_d)}\ls x\,\Big|\,\bp_d\right)\to P\left( \int_0^1w^2(t)B^2(t) \,dt \ls x\right)\quad a.s.
\end{align*}
for any continuous weight function $w(\cdot)$ with
\begin{align}
&	\lim_{t\to 0} t^{\alpha} w(t)<\infty,\quad \lim_{t\to 1}(1-t)^{\alpha}w(t)\quad\text{for some }0\ls \alpha <1/2,\notag\\
&\sup_{\eta\ls t\ls 1-\eta}w(t)<\infty\qquad\text{for all } 0<\eta\ls \frac 1 2.\label{eq_weights}
\end{align}
\end{cor}

\begin{figure}
\begin{center}
\subfloat[][Known variance]{\includegraphics[width=0.45\textwidth]{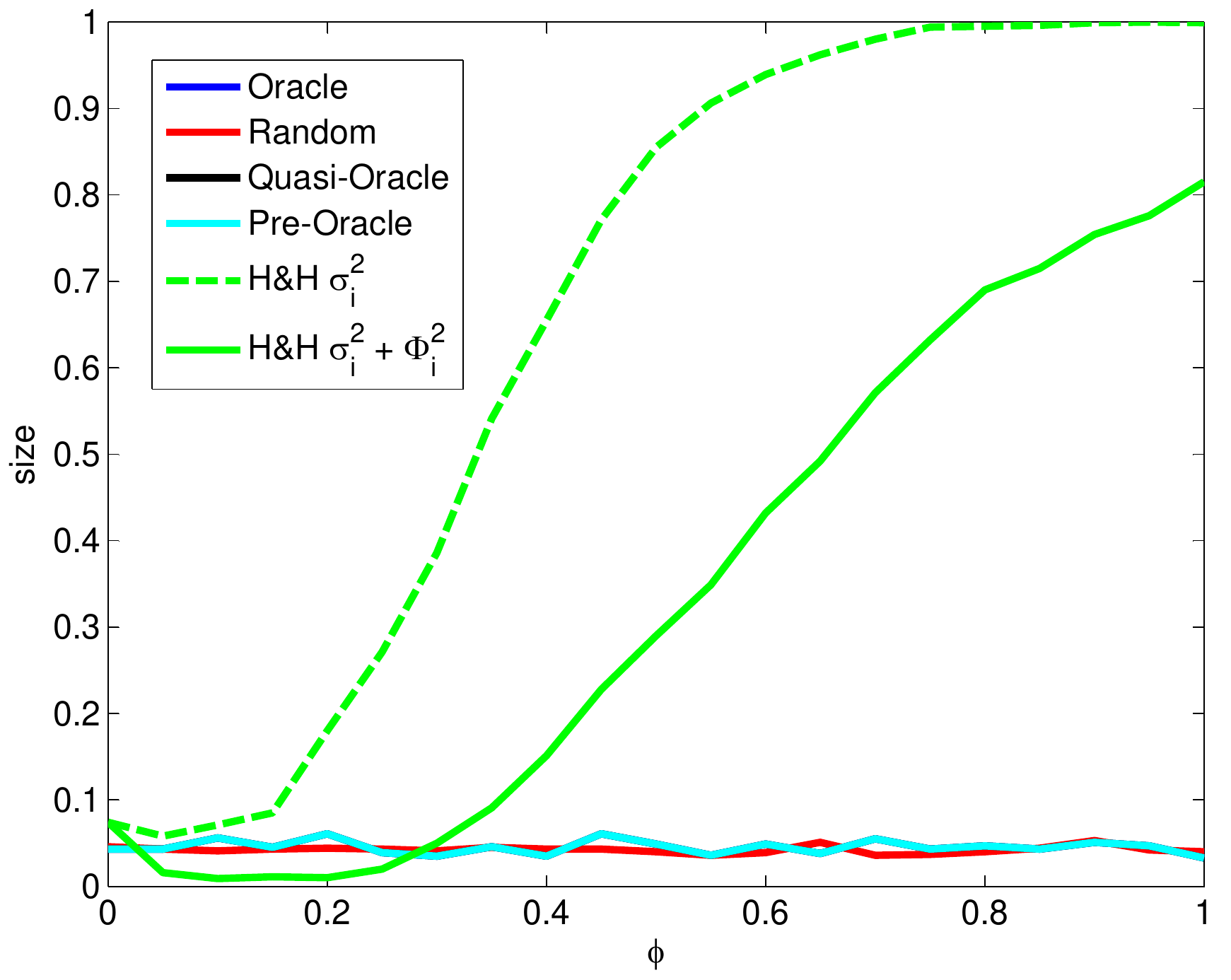}}\\
\subfloat[][Estimated variance as in \eqref{eq_est_tau}]{\includegraphics[width=0.45\textwidth]{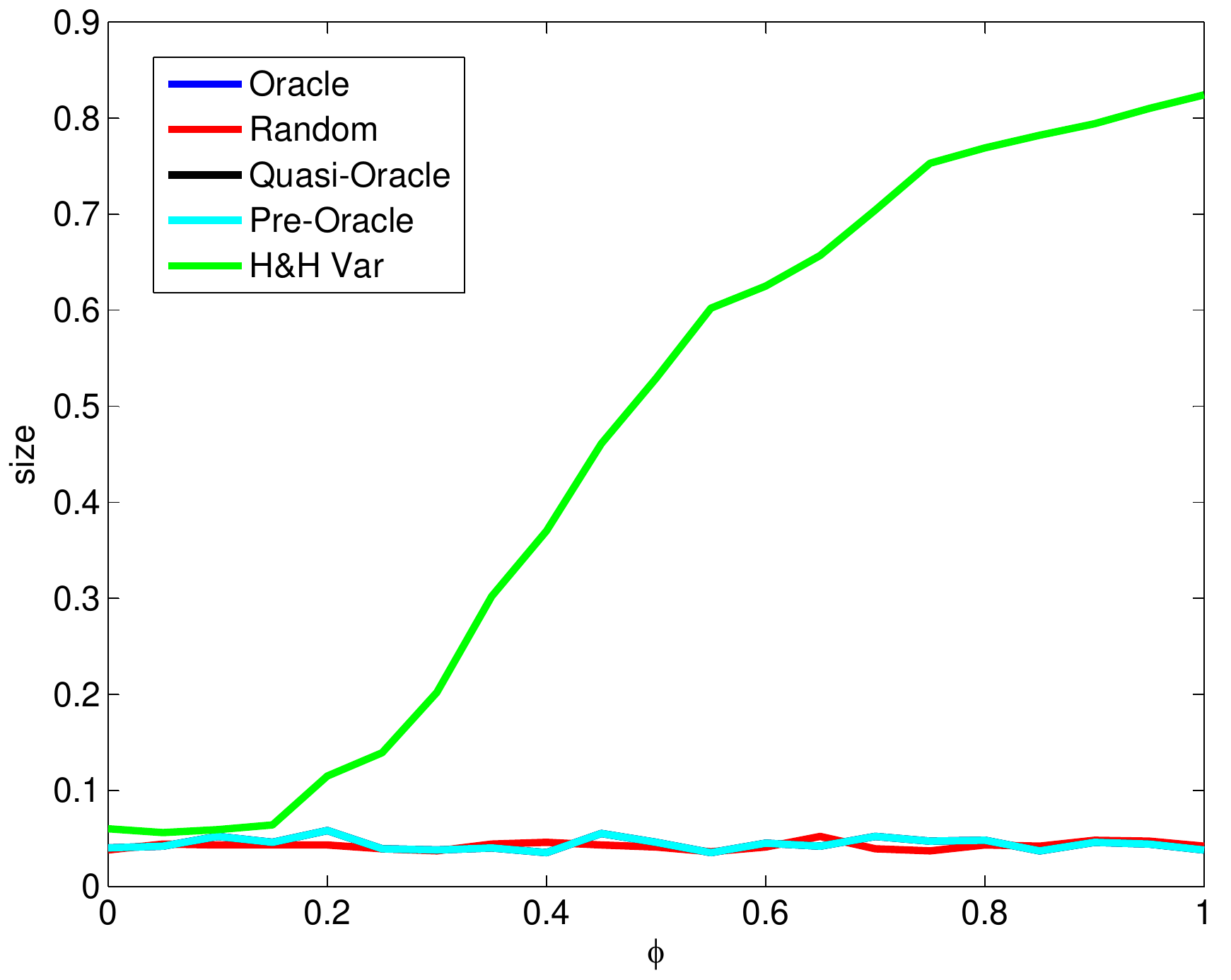}}
\subfloat[][Estimated variance as in \eqref{eq_est_tau2}]{\includegraphics[width=0.45\textwidth]{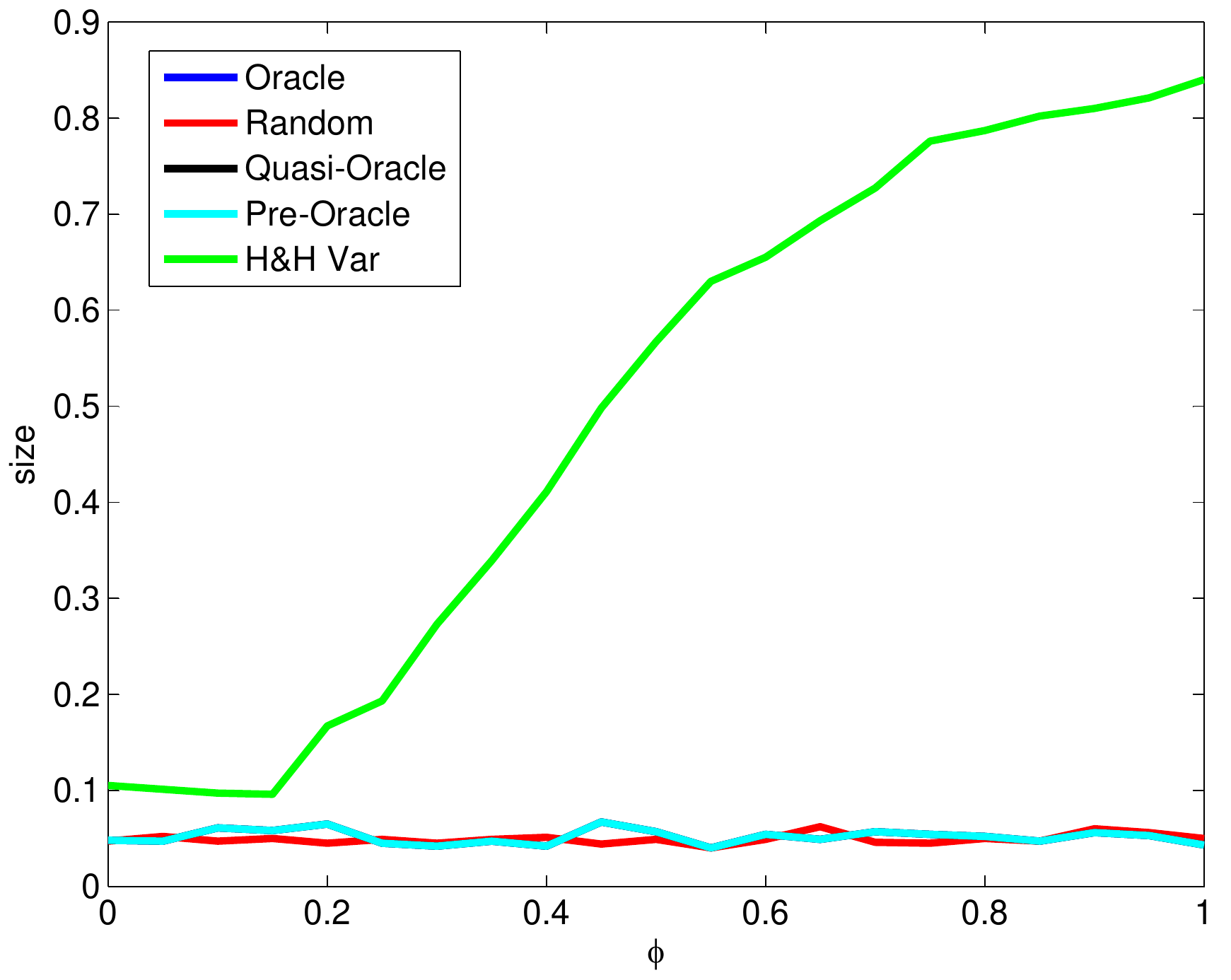}}\\
\caption{Size of tests as the degree of dependency between the components increases. As can be seen, all the projection methods, Oracle, Quasi-Oracle, Pre-Oracle and Random projections defined in Section \ref{section_oracle} maintain the size of the tests. Those based on using the full information as described in Section \ref{section_comparison} have size problems as the degree of dependency increases. The simulations correspond to Case $\mathcal{C}$.\ref{case_mixed} with $s_j=1, \Phi_j=\phi$, $j=1,\ldots,d$ 
with $d=200$, where $\phi$ is given on the x-axis).}\label{f:size}
\end{center}
\end{figure}

As can be seen in Figure \ref{f:size}, regardless of whether the variance is known or estimated, the projection methods all maintain the correct size even when there is a high degree of dependence between the different components (the specific projection methods and indeed the non-projection methods will be characterised in Section~\ref{section_oracle} below). The full tests, where size is not controlled, will be discussed in Section \ref{section_comparison}.

\subsubsection{Absolute high dimensional efficiency}

We are now ready to derive the high dimensional efficiency of projection statistics. Furthermore, we show that a related estimator for the location of the change is asymptotically consistent.

\begin{theorem}\label{th_ind_cont}
	Let the assumptions of Theorem~\ref{th_proj_null} either a) or b) on the errors respectively $\bp_d$ be fulfilled.	Furthermore, consider a weight function $w(\cdot)$ as in Corollary~\ref{th_asym_null_stat} fulfilling $w^2(x)\,\left( \int_0^xg(t)\,dt-x\int_0^1g(t)\,dt \right)^2\neq 0$. Then, both the max and sum statistic from Corollary~\ref{th_asym_null_stat} have the following absolute high dimensional efficiency:
	\begin{equation}\label{eq_22}
		\HDE_1(\bD_d,\bp_d):=	\frac{\|\boldsymbol{\Delta}_d\| \|\bp_d\| \cos(\alpha_{\boldsymbol{\Delta}_d,\bp_d})}{\tau(\bp_d)}= \frac{|\langle \bD_d,\bp_d\rangle|}{\tau(\bp_d)},
	\end{equation}
	where $\tau^2(\bp_d)$ is as in \eqref{eq_tau} and $\alpha_{\mathbf{u},\mathbf{v}}$ is the (smallest) angle between $\mathbf{u}$ and $\mathbf{v}$.
In addition, the asymptotic power increases with increasing multiplicative constant.
\end{theorem}

In the following, $\HDE_1(\bD_d,\bp_d)$ is fixed to the above representative of the class, so that different projection procedures with the same rate but with different constants can be compared.
\begin{rem}
	For random projections the high dimensional efficiency is a random variable. The convergences in Definition~\ref{def_HDE} is understood given the projection vector $\bp_d$, where we get either $a.s.$-convergence or $P$-stochastic convergence depending on whether  $\sqrt{T}\,\HDE_1(\bD_d,\bp_d)$ converges $a.s.$ or in a $P$-stochastic sense (in the latter case the assertion follows from the subsequence-principle).
\end{rem}

\begin{rem}
In particular all deviations from a stationary mean are detectable with asymptotic power one for weight functions $w>0$ as $\int_0^xg(t)\,dt-x\int_0^1g(t)\,dt\neq 0$ for non-constant $g$. It is this $g$ function which determines which weight function gives best power.

We derive the high dimensional efficiency for a given $g$ and disappearing magnitude of the change $\|\bD_d\|$. For an epidemic change situation with $g(x)=1_{\{\vth_1<x\ls \vth_2\}}$ for some $0<\vth_1<\vth_2<1$, this means that the duration of the change is relatively large but the magnitude relatively small with respect to the sample size. Alternatively, one could also consider the situation, where the duration gets smaller asymptotically~(see e.g.\ \cite{frick2014multiscale}) resulting in a different high dimensional efficiency, which is equal for both  the projection as well as the multivariate or panel statistic, as long as the same weight function and the same type of statistic (sum/max) is used.
Some preliminary investigations suggest that in this case using projections based on principle component analysis similar to \citet{AK1} can be advantageous, however this is not true for the setting discussed in this paper.
\end{rem}

The above result shows in particular that sufficiently large changes (as defined by the high dimensional efficiency) are detected asymptotically with power one. For such changes in the at-most-one-change situation, the following corollary shows that the corresponding change point estimator is consistent in rescaled time.

\begin{cor}\label{th_power_1_AMOC}
	Let the assumptions of Theorem~\ref{th_ind_cont} hold and additionally $\sqrt{T} \HDE_1(\bD_d,\bp_d)\to \infty$ $a.s.$
	Under the alternative of one abrupt change, i.e.\ $g(x)=1_{\{x>\vth\}}$ for some $0<\vth<1$, the estimator
	\begin{align*}
		\widehat{\vth}_T=\left\lfloor\frac{\arg\,max_k U^2_{d,T}(k/T)}{T}\right\rfloor
	\end{align*}
	is consistent for the change point in rescaled time, i.e.\
	\begin{align*}
		P\left(\left|\widehat{\vth}_T-\vth\right|\gs \eps\,|\,\bp_d\right) \to 0 \qquad a.s.
	\end{align*}
	An analogous statement holds, if the $\arg\max$ of $w^2(k/T)U^2_{d,T}(k/T)$ is used instead and $w^2(x)\left( (x-\vth)_+-x (1-\vth) \right)^2$ has a unique maximum at $\vth$, which is the case for many standard weight functions such as $w(t)=(t(1-t))^{-\beta}$ for some $0\ls \beta <1/2$.
\end{cor}

In the next section we will further investigate the high dimensional efficiency and see that the power depends essentially on the angle between $\Sigma^{1/2}\bp_d$ and the  'standardized' change $\Sigma^{-1/2}\bD$ if $\Sigma$ is invertible. In fact, the smaller the angle the larger the power. Some interesting insight will also come from the situation where $\Sigma$ is  not invertible by considering case $\mathcal{C}$.\ref{case_dep} above.

\subsection{High dimensional efficiency of oracle and random projections}\label{section_oracle}
In this section, we will further investigate the high dimensional efficiency of certain particular projections that can be viewed as benchmark projections.
In particular, we will see that the  efficiency depends only on the angle between the projection and the change both properly scaled with the underlying covariance structure.

The highest efficiency is obtained by $\mathbf{o}=\Sigma^{-1}\bD_d$ as the next theorem shows,  which will be called the oracle projection. This oracle is equivalent to a projection after first standardizing the data on the 'new' change $\Sigma^{-1/2}\bD_d$.  The corresponding test is related to the likelihood ratio statistic for i.i.d.\ normal innovations, where both the original mean and the direction (but not magnitude) of the change are known. As a lower (worst case) benchmark we consider
a scaled random projection $\br_{d,\Sigma}=\Sigma^{-1/2}\br_d$, where $\br_d$ is a random projection on the $d$-dimensional unit sphere. This is equivalent to a random projection onto the unit sphere after standardizing the data.
Both projections depend on $\Sigma$ which is usually not known so that it  needs to be estimated. The latter is rather problematic in particular in high dimensional settings without additional parametric or sparsity assumptions (see \citet{ZouHT2006}, \citet{BickelL2008} and \citet{FanLM2013} including related discussion). Furthermore, it is actually the inverse that needs to be estimated which results in additional numerical problems if $d$ is large. For this reason we check the robustness of the procedure with respect to not knowing { or misspecifying} $\Sigma$ in a second part of this section

In Section~\ref{section_comparison} we will  compare the efficiency of the above projections with a procedure taking the full information into account. We will show that we lose an order ${d}^{1/4}$ in terms of high dimensional efficiency between the oracle and the full panel data statistic and another ${d}^{ 1/4}$ between the panel and the random projection.

\subsubsection{Correctly scaled projections}

The following proposition characterizes which projection yields an optimal high dimensional efficiency associated with the highest power.

\begin{figure}
\begin{center}
\includegraphics[width=0.9\textwidth]{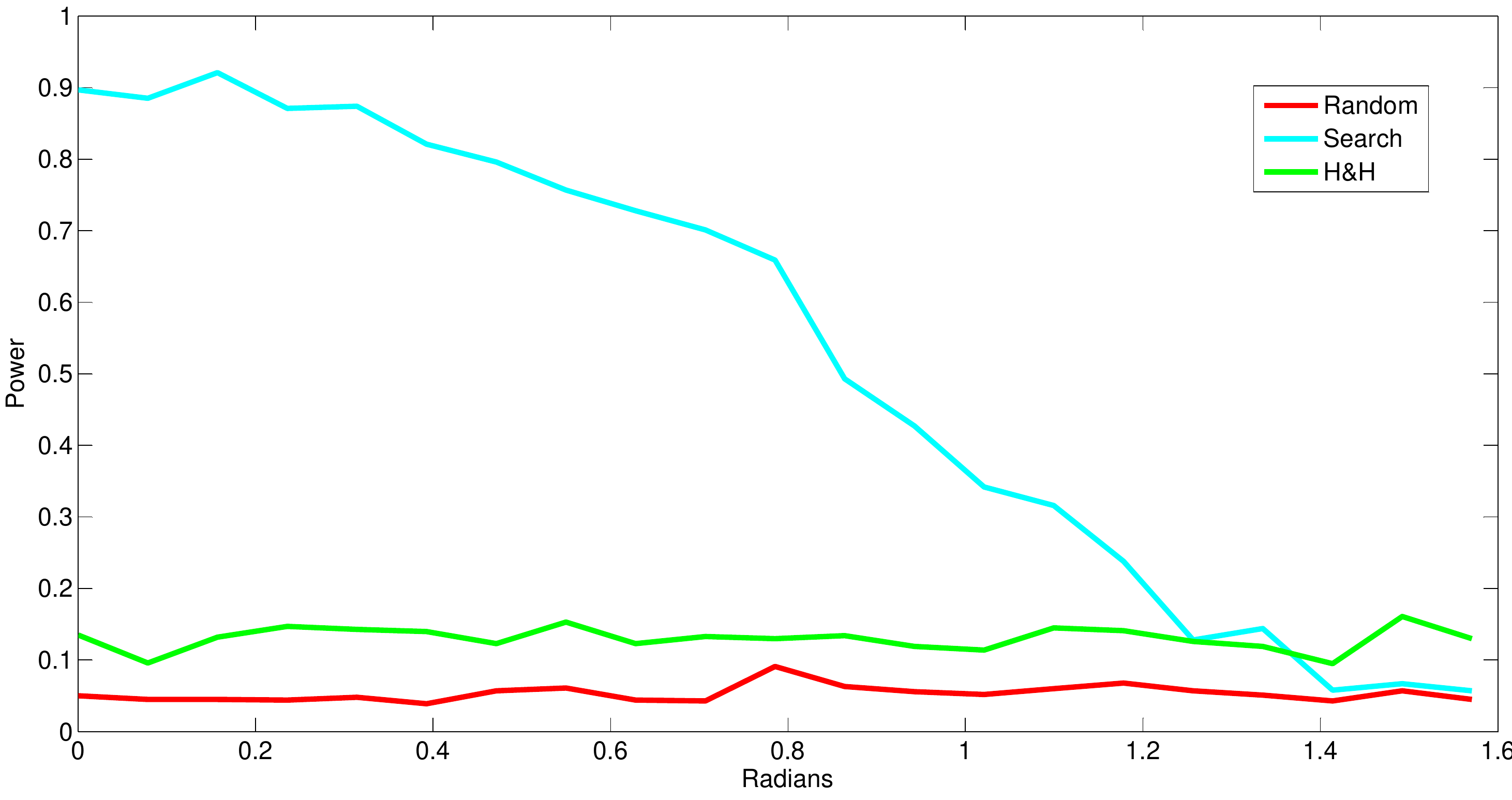}
\caption{Power of tests as the angle between the search direction and the oracle increases. As can be seen, the search projection method decreases similarly to cosine of the angle, while the random projection and \hor-\hus tests as introduced in Section \ref{section_comparison} are given for comparison. (Here $\Sigma_d = I_d$, $d=200$, and $\bD_d = 0.05\;\mathbf{1_d}$, corresponding to Case $\mathcal{C}$.\ref{case_ind}).}\label{f:angle}
\end{center}
\end{figure}

\begin{prop}\label{lemma_oracle_representation}
If $\Sigma$ is invertible, then
	\begin{equation}
		\HDE_1(\bD,\bp_d)=  \|\Sigma^{-1/2}\bD_d\|\cos(\alpha_{\Sigma^{-1/2}\bD_d,\Sigma^{1/2}\bp_d}).		\end{equation}
\end{prop}
Proposition~\ref{lemma_oracle_representation} shows in particular, that after standardizing the data, i.e.\ for $\Sigma=I_d$, the power depends solely on the cosine of the angle between the oracle and the projection (see Figure \ref{f:angle}).

From the representation in this proposition it follows immediately that the 'oracle' choice for the projection to maximize the high dimensional efficiency is  $\boldsymbol{o}=\Sigma^{-1}\boldsymbol{\Delta}_d$ as it maximizes the only term which involves the projection namely $\cos(\alpha_{\Sigma^{-1/2}\bD_d,\Sigma^{1/2}\bp_d})$.
Therefore, we define:
\begin{definition}\label{def_oracle}
		The projection $\bo=\Sigma^{-1}\boldsymbol{\Delta}_d$ is called {\bf oracle} if $\Sigma^{-1}$ exists.
	Since the projection procedure is invariant under multiplications with non-zero constants of the projected vector, all non-zero multiples of the oracle have the same properties, so that they correspond to a class of projections.
\end{definition}
By Proposition~\ref{lemma_oracle_representation} the oracle choice leads to a high dimensional efficiency of $\HDE_1(\bD_d,\bo)=\|\Sigma^{-1/2}\bD_d\|$.

Another way of understanding the Oracle projection is the following: If we first standardize the data, then for a projection on a unit (w.l.o.g.) vector the variance of the noise is constant and the signal is given by the scalar product of $\Sigma^{-1/2}\boldsymbol{\Delta}$ and the (unit) projection vector, which is obviously maximized by a projection with $\Sigma^{-1/2}\boldsymbol{\Delta}/\|\Sigma^{-1/2}\boldsymbol{\Delta}\|$ which is equivalent to using $\bp_d=\Sigma^{-1}\boldsymbol{\Delta}$ as a projection vector for the original non-standardized version.

So, if we know $\Sigma$ and want to maximize the efficiency respectively power close to a particular search direction $\bs_d$ of our interest, we should use the
{\bf scaled search direction} $\bs_{\Sigma,d}=\Sigma^{-1}\bs_d$ as a projection.

Because the cosine falls very slowly close to zero, the efficiency will be good if the search direction is not too far off the true change. From this, one could get the impression that even a scaled random projection $\br_{\Sigma,d}=\Sigma^{-1/2}\br_d$  may not do too badly, where $\br_d$ is a uniform random projection on the unit sphere.
This is equivalent to using a random projection on the unit sphere after standardizing the data, which also explains the different scaling as compared to the oracle or the scaled search direction, where the change $\bD_d$ is also transformed to $\Sigma^{-1/2}\bD_d$ by the standardization.
However, since for increasing $d$ the space covered by the far away angles is also increasing, the high dimensional efficiency of the scaled random projection is not only worse than the oracle by a factor $\sqrt{d}$ but also by a factor ${d^{1/4}}$  than a full multivariate or panel statistic which will be investigated in detail in Section~\ref{section_comparison}.

The following theorem shows the high dimensional efficiency of the scaled random projection.

\begin{theorem}\label{th_random_proj_scaled}
	Let the alternative hold, i.e.\ $\|\bD_d\|\neq 0$. Let $\br_d$ be a random uniform projection on the $d$-dimensional unit sphere and $\br_{\Sigma,d}=\Sigma^{-1/2}\br_d$, then for all $\epsilon>0$ there exist  constants $c,C>0$, such that
\begin{align*}
	P\left(c\ls   \HDE^2_1(\bD_d,\br_{\Sigma,d})\,\frac{d}{\|\Sigma^{-1/2}\,\bD_d\|^2} \ls C\right)\gs 1-\epsilon.
\end{align*}
\end{theorem}

Such a random projection on the unit sphere can be obtained as follows: Let $X_1,\ldots,X_d$ be i.i.d.\ N(0,1), then $\br_d=(X_1,\ldots,X_d)^T/\|(X_1,\ldots,X_d)^T\|$ is uniform on the $d$-dimensional unit sphere \cite{Marsaglia72}.

Comparing the high dimensional efficiency of the scaled random projection with the one obtained for the oracle projection (confer Proposition~\ref{lemma_oracle_representation}) it becomes apparent that we lose an order $\sqrt{d}$. In Section~\ref{section_comparison} we will see that the panel statistic taking the full multivariate information into account has a contiguous rate just between those two losing $d^{1/4}$ in comparison to the oracle but gaining ${d^{1/4}}$ in comparison to a scaled random projection.
From these results one obtains a cone around the search direction such that the projection statistic has higher power than the panel statistic, if the true change falls within this cone.

\begin{figure}
\begin{center}
\includegraphics[width=0.9\textwidth]{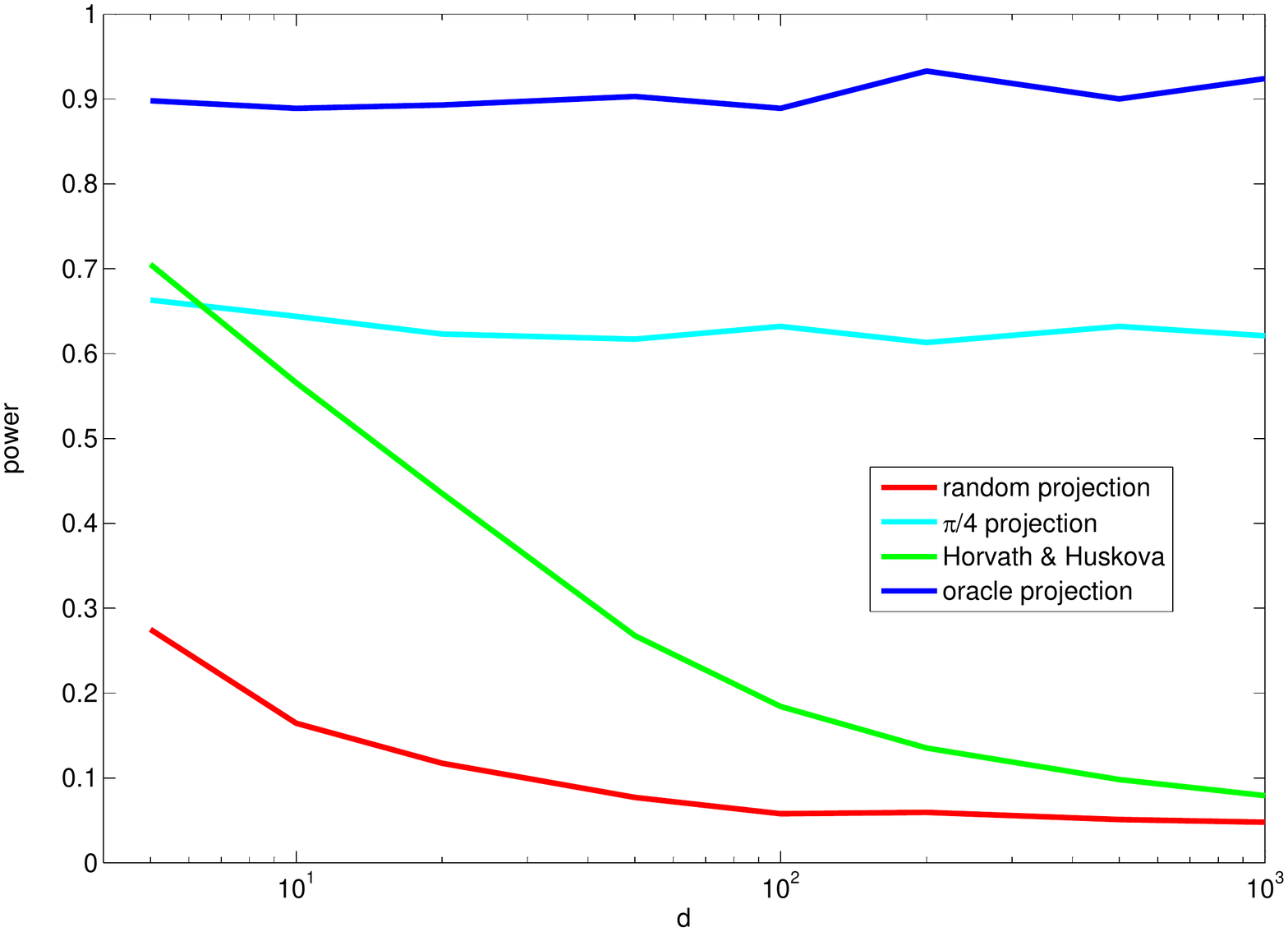}
\caption{Power of the tests as $d$ increases with a fixed sample size ($T=100$). Here $ \|\bD_d\| = \mbox{const.}$
and $\Sigma_d = I_d$, i.e.\ $\|\Sigma^{-1/2}\bD_d\|= \mbox{const.}$, corresponding to Case $\mathcal{C}$.\ref{case_ind}. This gives roughly constant power for fixed angle projection tests (as $\|\bD_d\|$ is constant), while results in decreasing power for both the panel statistic test and random projections as predicted by theory.}\label{f:d_increasing}
\end{center}
\end{figure}

Figure \ref{f:d_increasing} shows the results of some simulations showing that a change that can be detected for the oracle with constant power as $d$ increases rapidly loses power for the panel statistic as  predicted by its high dimensional efficiency in Section \ref{section_comparison} as well as for the random projection. This and the following simulations show clearly that the concept of high dimensional efficiency is indeed capable of explaining the small sample power of a statistic very well.

\subsubsection{The oracle in the case of non-invertibility}\label{s:noninvert}

Let us now  have a look at the situation if $\Sigma$ is not invertible hence the above oracle does not exist. To this end, let us consider Case $\mathcal{C}$.2 above -- other non-invertible dependent situations can essentially be viewed in a very similar fashion, but become a combination of the two scenarios below.

\begin{figure}
\begin{center}
\subfloat[][Angle between $\bD_d$ and $\Phi$ = 0 radians]{\includegraphics[width=0.45\textwidth]{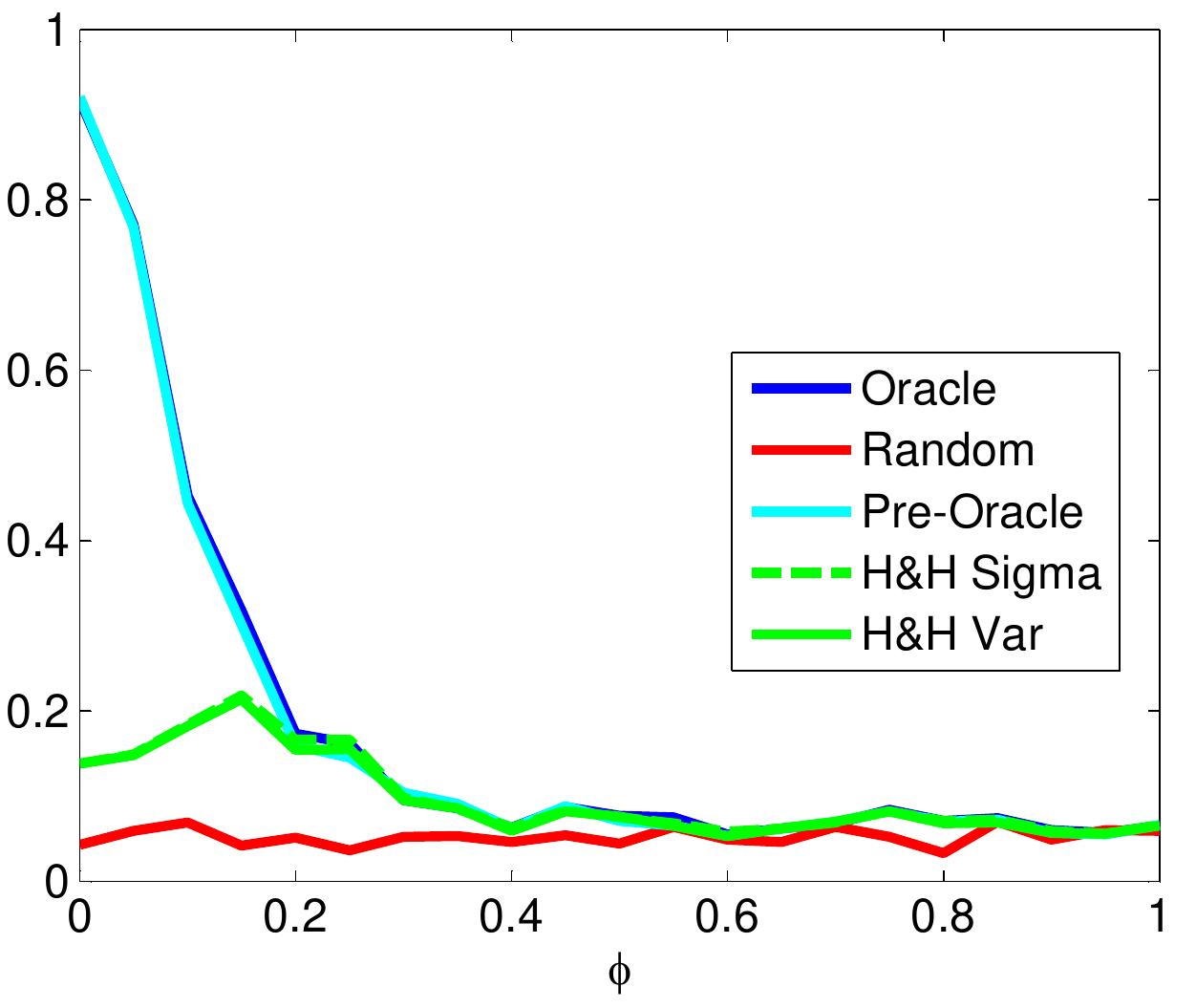}}
\subfloat[][Angle between $\bD_d$ and $\Phi$ = $\pi$/8 radians]{\includegraphics[width=0.45\textwidth]{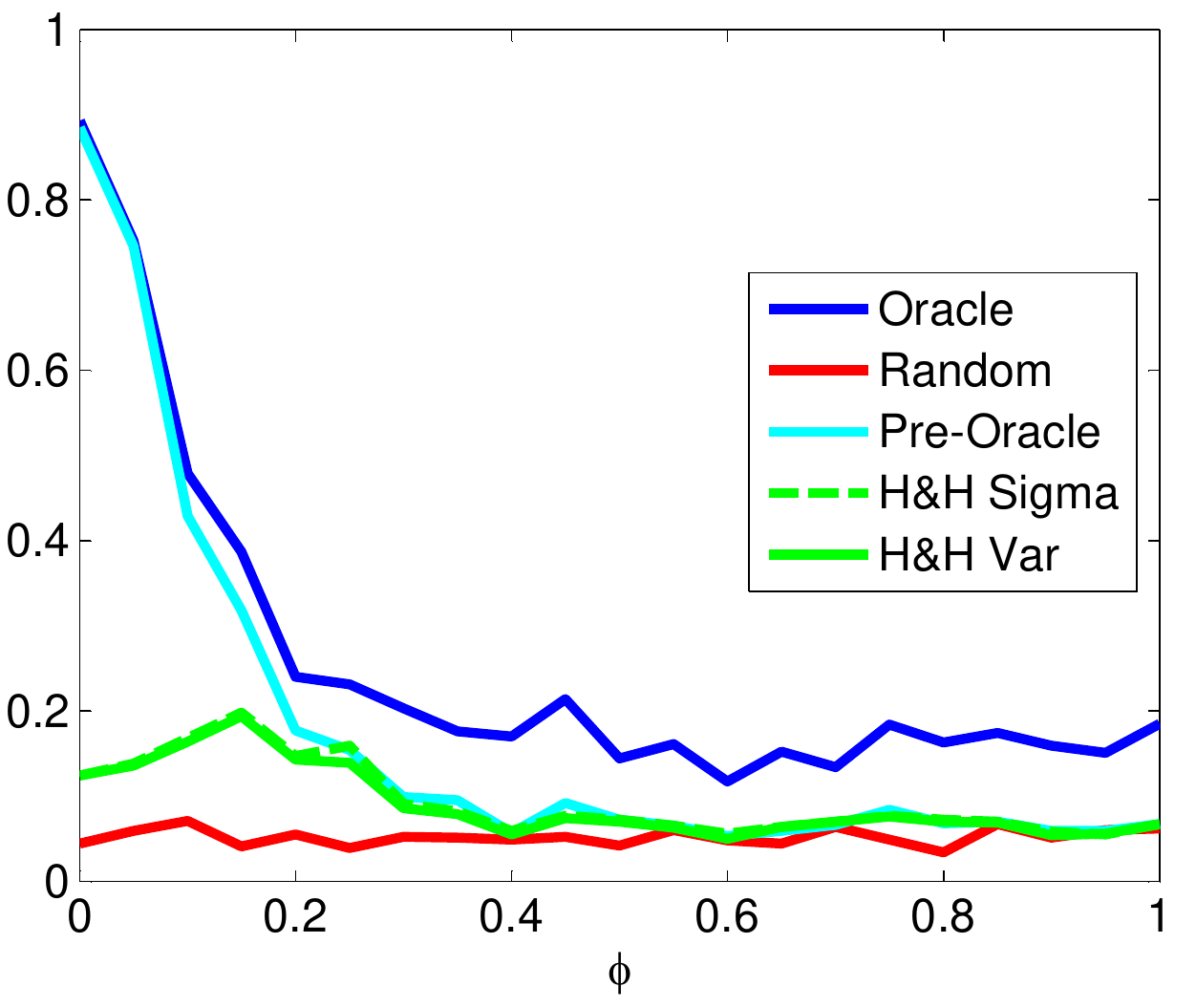}}\\
\subfloat[][Angle between $\bD_d$ and $\Phi$ = $\pi$/4 radians]{\includegraphics[width=0.45\textwidth]{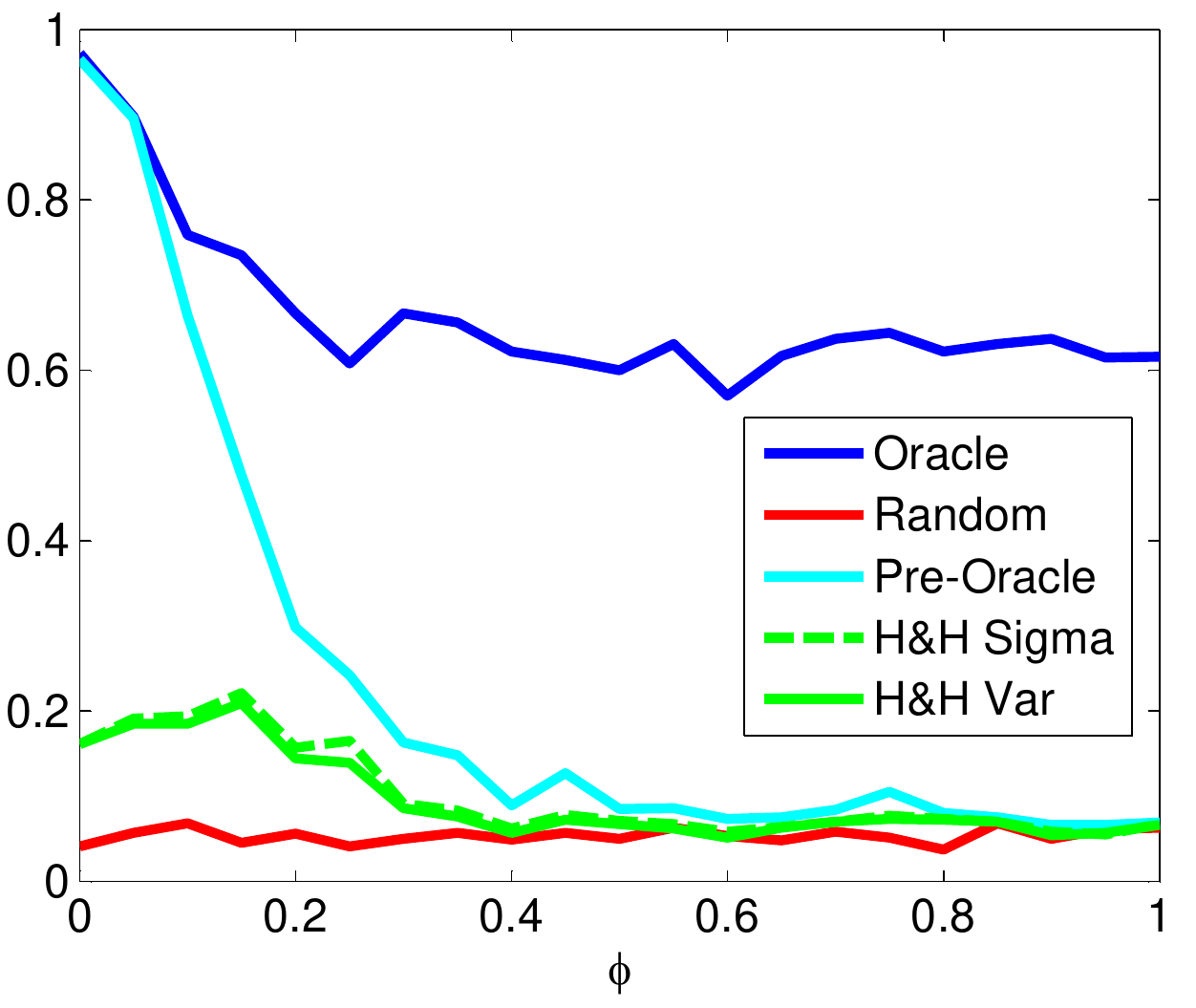}}
\subfloat[][Angle between $\bD_d$ and $\Phi$ = $\pi$/2 radians]{\includegraphics[width=0.45\textwidth]{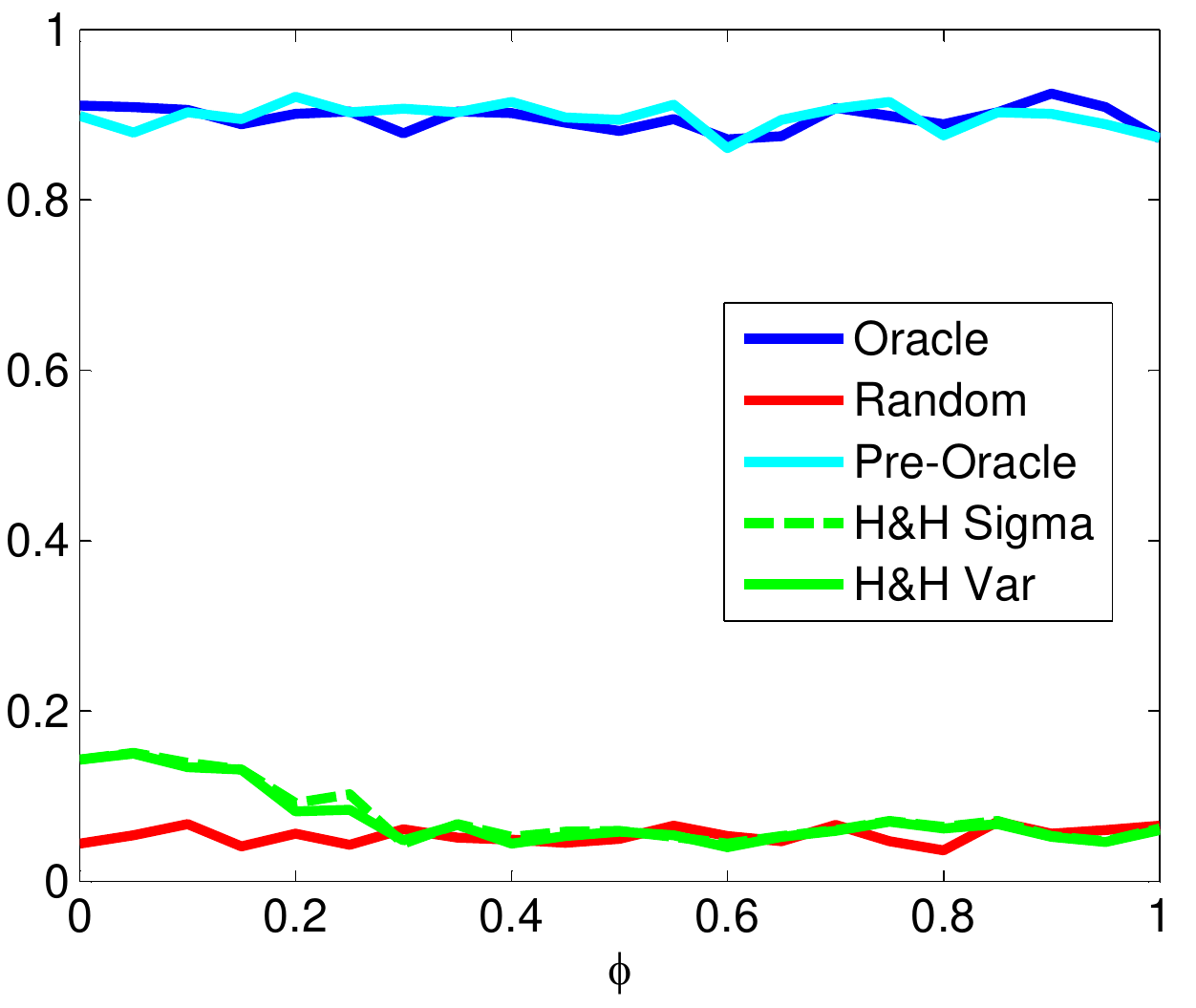}}\\
\caption{Power of tests as the angle between the change and the direction of dependency increases. As can be seen, if the change lies in the direction of dependency, then all methods struggle, which is in line with the theory of Section \ref{section_oracle}. However, if the change is orthogonal to the dependency structure the projection method works regardless of whether the dependency is taken into account or not.  H\&H Sigma and Var as in Section \ref{section_comparison} represent the panel tests taking into account the true or estimated variances of the components. All results are empirically size corrected to account for the size issues seen in Figure \ref{f:size}. 
($s_j=1$, $\Phi_j=\phi$, $j=1,\ldots,d$ with $d=200$, $\|\bD_d\|=0.05\sqrt{d}$, corresponding to Case $\mathcal{C}$.\ref{case_mixed}), with $\phi$ as given on the x-axis.}\label{f:HHcomp2}
\end{center}
\end{figure}

\setcounter{case}{1}
\begin{case}[Fully dependent Components]
	In this case $\Sigma=\bPhi_d\bPhi_d^T$ is a rank 1 matrix and not invertible. Consequently, the oracle as in Definition~\ref{def_oracle} does not exist. To understand the situation better, we have to distinguish two scenarios:
			\begin{enumerate}[(i)]
				\item If $\bPhi_d$ is not a multiple of $\bD_d$ we can transform the data into a noise-free sequence that only contains the signal by projecting onto a vector that is orthogonal to $\bPhi_d$ (cancelling the noise term) but not to $\bD_d$. All such projections are in principle equivalent as they yield the same signal except for a different scaling which is not important if there is no noise present. Consequently, all such transformations could be called oracle projections.
					\item	On the other hand if $\bD_d$ is a multiple of $\bPhi_d$, then any projection cancelling the noise will also cancel the signal. Projections that are orthogonal to $\bPhi_d$ hence by definition also to $\bD_d$ will lead to a constant deterministic sequence hence to a degenerate situation.  All other projections lead to the same (non-degenerate) time series except for multiplicative constants and different means (under which the proposed change point statistics are invariant by definition) so all of them could be called oracles.
	\end{enumerate}
	The following interpretation also explains the above mathematical findings: In this situation, all components are obtained from one common factor $\{\eta_t\}$ with different weights according to $\bPhi_d$ i.e.\ they move in sync with those weights. If a change is proportional to $\bPhi_d$ it could either be attributed to the noise coming from $\{\eta_t\}$ or from a change, so it will be difficult to detect as we are essentially back in a duplicated rank one situation and no additional information about the change can be obtained from the multivariate situation. However, if it is not proportional to $\bPhi$, then it is immediately clear (with probability one) that a change in mean must have occurred (as the underlying time series no longer moves in sync). This can be seen to some extent in Figure \ref{f:HHcomp2}, where the different panels in the figure mimic the different scenarios as outlined above (with a large value of $\phi$ being close to the non-invertible situation).
\end{case}

\subsubsection{Misscaled projections with respect to the covariance structure}\label{section_misscaled_oracle}

The analysis in the previous section requires the knowledge or a precise estimate of the inverse (structure) of $\Sigma$. However, in many situations such an estimate may not be feasible or too imprecise due to one or several of the below reasons, where the problems get worse due to the necessity for inversion \begin{itemize}
	\item If $d$ is large in comparison to $T$ statistical estimation errors can accumulate and identification may not even be possible \cite{BickelL2008}.
	\item The theory can be generalized to time series errors but in this case the covariance matrix has to be replaced by the long-run covariance (which is proportional to the spectrum at 0) and is much more difficult to estimate \cite{AK2,KirchTK2012}.
	\item Standard covariance estimators will be inconsistent under alternatives as they are contaminated by the change points. Consequently, possible changes have to be taken into account, but even in a simple at most one change situation it is unclear how best to generalize the standard univariate approach as in \eqref{eq_est_tau2} as opposed to \eqref{eq_est_tau} to a multivariate situation as the estimation of a joint location already requires an initial weighting for the projection (or the multivariate statistic). Alternatively, component-wise univariate estimation of the change points could be done but require a careful asymptotic analysis in particular in a setting with $d\to\infty$.
	\item If $d$ is large, additional numerical errors may arise when inverting the matrix \cite[Ch 14]{Higham2002}.
\end{itemize}

We will now investigate the influence of misspecification or estimation errors on the high dimensional efficiency of a {\bf misscaled oracle} $\bo_M=\bM^{-1} \bD_d$ in comparison to the {\bf misscaled random projection} $\br_{\bM,d}=\bM^{-1/2}\br_d$, where we only assume that the assumed covariance structure $\bM$ is symmetric and positive definite and model $\mathcal{A}$.\ref{model_VAR0} is fulfilled.

\begin{theorem}\label{th_random_proj}
	Let the alternative hold, i.e.\ $\|\bD_d\|\neq 0$. Let $\br_d$ be a random projection on the $d$-dimensional unit sphere and $\br_{\bM,d}=\bM^{-1/2}\br_d$ be the misscaled random projection. Then, there exist for all $\epsilon>0$ constants $c,C>0$, such that
	\begin{align*}
		&P\left(c\ls  \HDE^2_1(\bD_d,\br_{\bM,d})\,\frac{\mbox{tr}\left(\bM^{-1/2}\Sigma \bM^{-1/2}\right) }{\|\bM^{-1/2}\bD_d\|^2}\ls C  \right)\gs 1-\epsilon,
	\end{align*}
	where $\mbox{tr}$ denotes the trace.\end{theorem}

We are now ready to prove the main result of this section stating that the high dimensional efficiency of a misscaled oracle can never be worse than the corresponding misscaled random projection.

\begin{theorem}\label{theorem_oracle_rand}
Let Assumption~$\mathcal{A}$.\ref{model_VAR0} hold.
Denote  the misscaled oracle by $\bo_M=\bM^{-1} \bD_d$, then
\begin{align*}
	\HDE_1^2(\bD_d, \bo_M)\gs \frac{\|\bM^{-1/2}\bD_d\|^2}{\mbox{tr}(\bM^{-1/2}\Sigma\bM^{-1/2})}
		\end{align*}
		where $\mbox{tr}$ denotes the trace and equality holds iff there is only one common factor which is weighted proportional to  $\bD_d$,
\end{theorem}

Because it is often assumed that components are independent and it is usually feasible to estimate the variances of each component,
we  consider the correspondingly misscaled oracles, which are scaled with the identity matrix (pre-oracle) respectively with the diagonal matrix of variances (quasi-oracle). The quasi-oracle is of particular importance as it uses the same type of misspecification as the panel statistic discussed in Section \ref{section_comparison} below.
\begin{definition}
	\begin{enumerate}[(i)]
		\item The projection  $\po=\bD_d$  is called {\bf pre-oracle}.
		\item The projection $\qo=\Lambda_d^{-1}\bD_d=(\delta_1/\sigma_1^2,\ldots,\delta_d/\sigma_d^2)^T$, $\Lambda_d=\mbox{diag}(\sigma_1^2,\ldots,\sigma_d^2)$  is called {\bf quasi-oracle}, if $\sigma_j^2>0$, $j=1,\ldots,d$.
		\end{enumerate}
As with the oracle, these projections should be seen as representatives of a class of projections.
	\end{definition}

	The following proposition shows that in the important special case of uncorrelated components, the (quasi-)oracle and pre-oracle have an efficiency of same order if the variances in all components are bounded and bounded away from zero. The latter assumption is also needed for the panel statistic below and means that all components are on similar scales. In addition, the efficiency of the quasi-oracle is even in the misspecified situation always better than an unscaled random projection.
\begin{prop}\label{prop_oracle}Assume that all variances are on the same scale, i.e.  there exist $c,C$ such that $0<c\ls \sigma_i^2<C<\infty$ for $i=1,\ldots,d$.
	\begin{enumerate}[a)]
		\item Let $\Sigma=\mbox{diag}(\sigma_1^2,\ldots,\sigma_d^2)$, then
	\begin{align*}
		& \frac{c^2}{C^2}\HDE^2_1(\bD_d,\qo)\ls \HDE^2_1(\bD_d,\po)\ls \HDE^2_1(\bD,\qo)=\|\Sigma^{-1/2}\bD_d\|^2.
	\end{align*}
\item Under Assumption~$\mathcal{A}$.\ref{model_VAR0}, it holds
	\begin{align*}
			\HDE_1^2(\bD_d,\qo)\gs \frac{c^2}{C^2}\,\frac{\|\bD_d\|^2}{\mbox{tr}(\Sigma)}.
		\end{align*}
	\end{enumerate}
		\end{prop}

			Figures \ref{f:HHcomp2} and \ref{f:degeneracy} show the results of some small sample simulations of case $\mathcal{C}.3$ confirming again that the theoretical findings based on the high dimensional efficiency are indeed able to predict the small sample power for the test statistics. Essentially, the following assertions are confirmed:
			\begin{enumerate}[1)]
				\item The power of the pre- and quasi-oracle is always better than the power of the misscaled random projection (the random projection assumes an identity covariance structure).
				\item The power of the  (correctly scaled) oracle can become as bad as the power of the (misscaled) random projection but only if $\bPhi_d\sim\bD_d$. In this case the power of the misscaled panel statistic  (i.e.\ where the statistic but not the critical values are constructed under the wrong assumption of independence between components) is equally bad.
				\item While the power of the (misscaled) panel statistic becomes as bad as the power of the (misscaled) random projection for $\phi\to\infty$ irrespective of the angle between $\bD_d$ and $\bPhi_d$, it can be significantly better for the pre- and quasi-oracle.
In fact, we will see in Section~\ref{section_comparison} that the high dimensional efficiency of the misspecified panel statistic will be of the same order as a random projection for any choice $\bPhi_d$ with $\bPhi_d^T\bPhi_d\sim d$, {irrespective of the direction of any change that might be present}.

			\end{enumerate}

We will now have a closer look at the three standard examples in order to understand the behavior in the simulations better (Case $\mathcal{C}.1$ is included in the simulations for $\Phi=0$, while $\mathcal{C}.3$ is the limiting case for $\Phi\to\infty$).

	\setcounter{case}{0}
\begin{case}[Independent components]
	If the components are uncorrelated, each with variance $\sigma_i^2$, i.e.\ $\Sigma_1=\mbox{diag}(\sigma_1^2,\ldots,\sigma_d^2)$, we get
			\begin{align*}
				\mbox{tr}(\Sigma_1)=\sum_{j=1}^d\sigma_j^2,
			\end{align*}
			which is of order $d$ if $0<c\ls \sigma_j^2\ls C<\infty$.
			Proposition~\ref{prop_oracle}, Theorem~\ref{th_random_proj_scaled} and Theorem~\ref{th_random_proj} show that in this situation both the high dimensional efficiency of the pre- and (quasi-)oracle are of an order $\sqrt{d}$ better than the correctly scaled and unscaled random projection.
		\end{case}
		The second case shows that high dimensional efficiency of misscaled oracles can indeed become as bad as a random projection and helps in the understanding of the mixed case:
		\begin{case}[Fully dependent components]
			As already noted in \ref{s:noninvert} we have to distinguish two cases:
			\begin{enumerate}[(i)]
			\item 	
				If $\bD_d$ is not a  multiple of $\bPhi_d$, then the power depends on the angle of the projection with $\bPhi_d$ with maximal power for an orthogonal projection. So the goodness of the oracles depends on their angle with the vector $\bPhi_d$.
			\item If $\bD_d$ is a multiple of $\bPhi_d$,  the pre- and quasi-oracle are not orthogonal to the change, hence they share the same high dimensional efficiency with any scaled random projection as all random projections are not orthogonal to $\bPhi_d$ with probability 1.
	\end{enumerate}
		\end{case}

We can now turn to the mixed case that is also used in the simulations.
		
		\begin{case}[Mixed case]
				Let $\ba_j=(0,\ldots,s_j,\ldots,0)^T$ the vector which is $s_j>0$ at point $j$ and zero everywhere else, and $\ba_{d+1}=\bPhi_d=(\Phi_1,\ldots,\Phi_d)^T$, $\ba_j=\boldsymbol 0$ for $j\gs d+2$.  Then $\Sigma_3=\mbox{diag}(s_1^2,\ldots,s_d^2)+\bPhi_d\bPhi_d^T$ and
				\begin{align}\label{eq_trace_mixed_case}
				\mbox{tr}(\Sigma_3)=\sum_{j=1}^ds_j^2+\sum_{j=1}^d\Phi_j^2.
			\end{align}

			The high dimensional efficiency of the pre-oracle can become as bad as for the random projection if the change $\bD_d$ is a multiple of the common factor $\bPhi_d$ and there is a substantial common effect. This is similar to Case~$\mathcal{C}$.\ref{case_dep} (which can be seen as a limiting case for increasing $\|\bPhi_d\|$).
			Intuitively, the problem is the following: By projecting onto the change, we want to maximize the signal { i.e.\ the change in the projected sequence} while minimizing the noise. In this situation however, the common factor dominates the noise in the projection as it essentially adds up in a linear manner, while the uncorrelated components add up only in the order of $\sqrt{d}$ (CLT). Now, projecting onto $\bD_d=\bPhi_d$ maximizes not only the signal but also the noise, which is why we cannot gain anything (but this also holds true for competing procedures { as in Section \ref{section_comparison} below}).
		
			More precisely, in $\mathcal{C}$.3 it holds $\tau^2(\po)=\sum_{j=1}^ds_j^2\delta_j^2+\left( \sum_{j=1}^d\delta_j\Phi_j \right)^2$. If additionally $\bD_d=k\bPhi_d$,  for some $k>0$, we get the following high dimensional efficiency for the pre-oracle	 by \eqref{eq_22}
			\begin{align*}
				\HDE_1(\bD_d,\po)=\frac{\|\bD_d\|}{\sqrt{\sum_{i=1}^ds_i^2\left(\frac{\delta_i}{\|\bD_d\|}\right)^2+\|\bPhi_d\|^2}}.
	\end{align*}
	The  high dimensional efficiency for the unscaled random projection is given by (confer Theorem~\ref{th_random_proj} and \eqref{eq_trace_mixed_case})
	\begin{align*}
		\frac{\|\bD_d\|}{\sqrt{\sum_{j=1}^ds_j^2+\|\bPhi_d\|^2}}.
	\end{align*}
As soon as $s_j, \Phi_j$ are of the same order, i.e.\ $0<c\ls s_j,\Phi_j\ls C<\infty$ for all $j$, the pre-oracle behaves as badly as	 the unscaled random projection. The same holds for the quasi-oracle under the same { assumptions.}
Interestingly, however, in this particular situation, even the oracle is of the same order as the random projection if the $s_j$ are of the same order, i.e.\ $0<c\ls s_j<C<\infty$. More precisely we get (for a proof we refer to the Section~\ref{section_proofs})
	\begin{align}\label{ex_oracle}
		 &\HDE_1(\bD_d,\bo)=\frac{\|\bD_d\|}{\sqrt{1+\sum_{j=1}^d\frac{\Phi_j^2}{s_j^2}}}\,\sqrt{\frac{\sum_{j=1}^d\frac{\delta_j^2}{s_j^2}}{\sum_{j=1}^d\delta_j^2}}.
	\end{align}
Figure \ref{f:HHcomp2} shows simulations which confirm the underlying theory in finite samples.

On the other hand, if $\bD_d$ is orthogonal to $\bPhi_d$, then the noise from $\bPhi_d$ cancels for the pre-oracle projection and we get the rate
	\begin{align*}
		\HDE_1(\bD_d,\po)=\frac{\|\bD_d\|}{\sqrt{\sum_{i=1}^ds_i^2\left(\frac{\delta_i}{\|\bD_d\|}\right)^2}},
	\end{align*}
	which is of the order $\|\bD_d\|^2$ if the $s_j$ are all of the same order. Anything between those two cases is possible and depends on the angle between $\bD$ and $\bPhi_d$ (again see Figures \ref{f:HHcomp2} and \ref{f:degeneracy} for finite sample simulations).

\begin{figure}
\begin{center}
\subfloat[][No Dependency - $\phi=0$]{\includegraphics[width=0.45\textwidth]{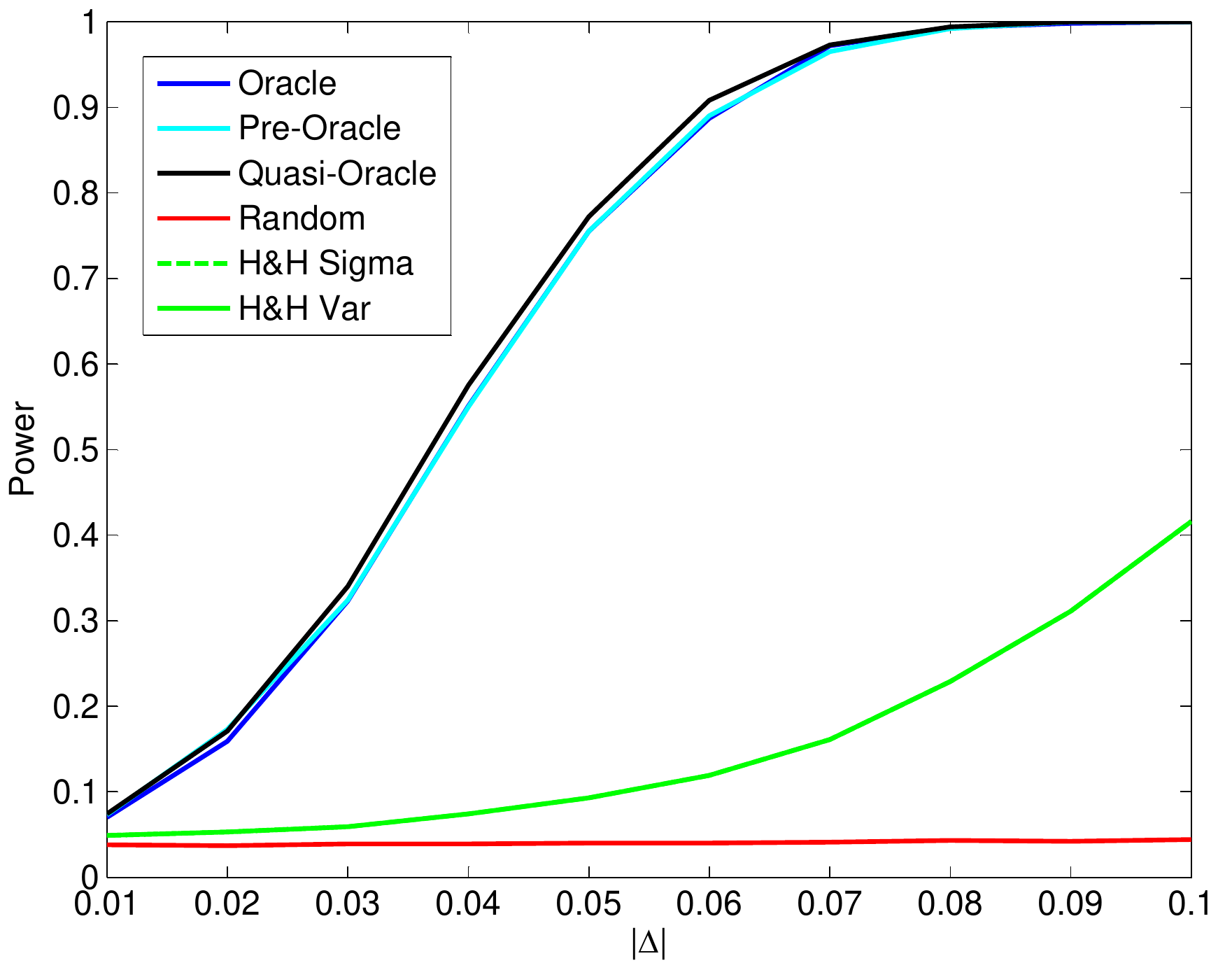}}
\subfloat[][$\phi=0.5$]{\includegraphics[width=0.45\textwidth]{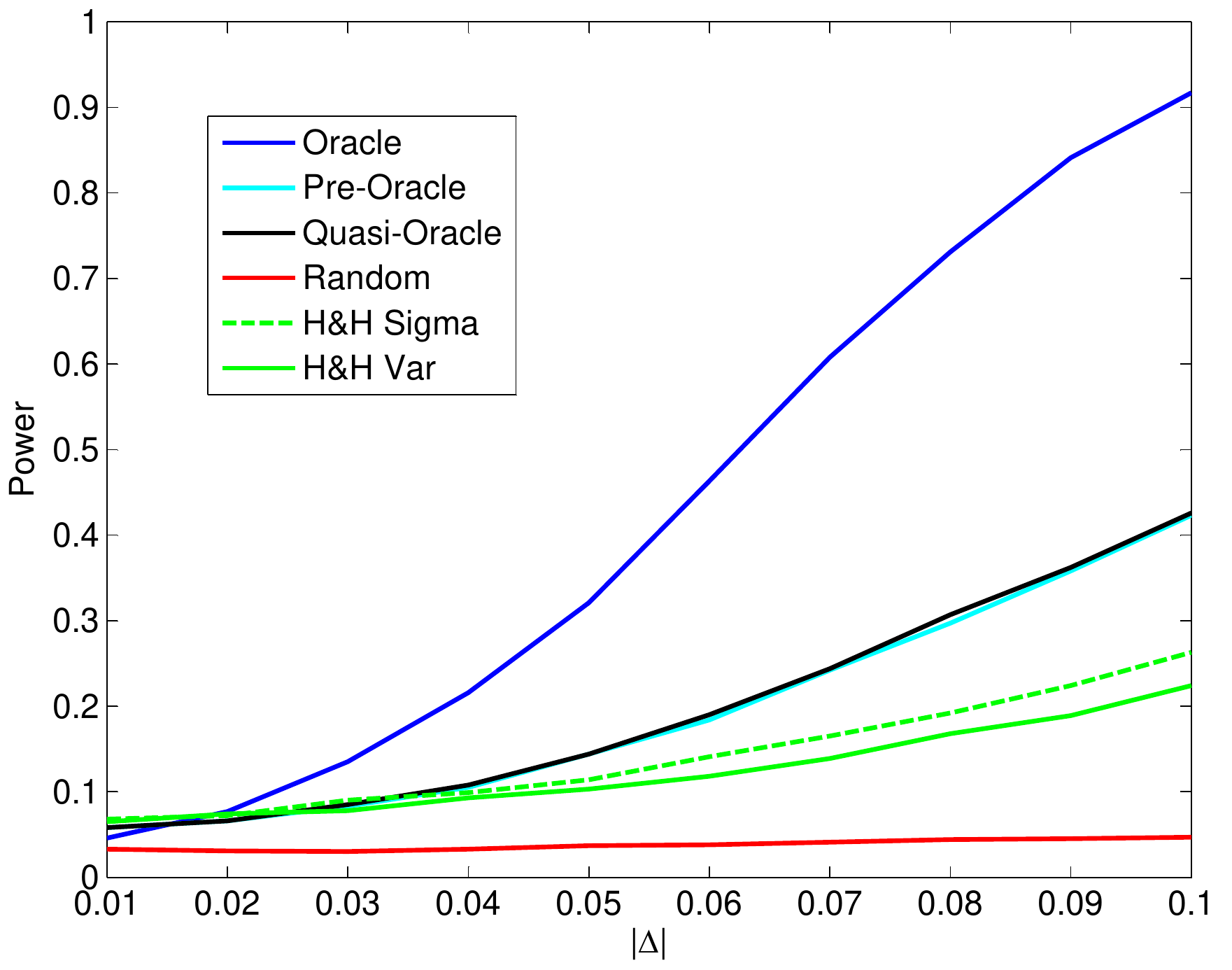}}\\
\subfloat[][$\phi=1$]{\includegraphics[width=0.45\textwidth]{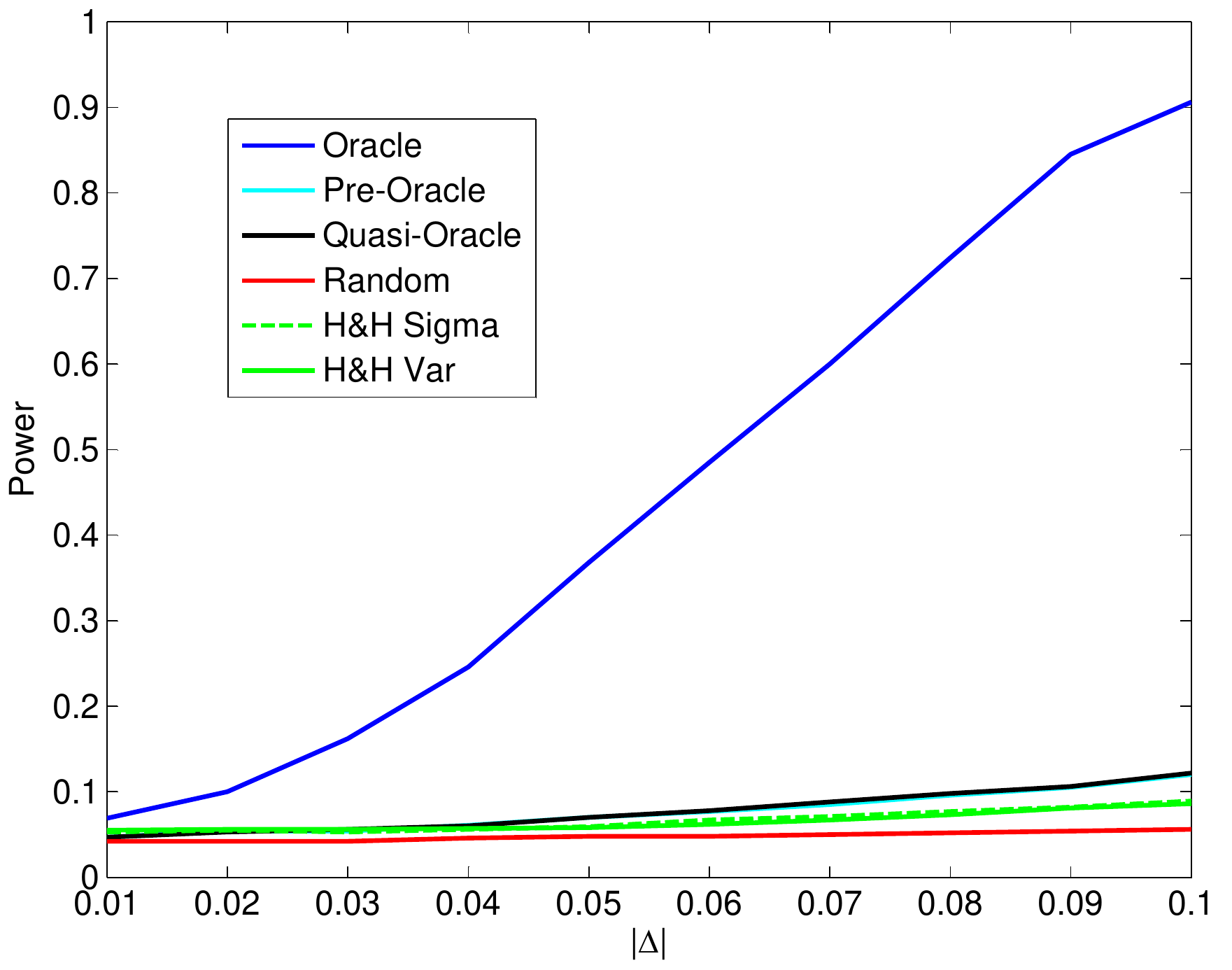}}
\caption{Power of tests as the dependency increases. The covariance structure becomes closer to degenerate across the three graphs, but in all cases the pre-oracle and quasi-oracle still outperform random projections, although they become closer as the degeneracy increases. Here different variances are used across components, namely $s_i=0.5+i/d$, $\Phi_i=\phi_i$, $i=1,\ldots,d$, 
$d=200$, angle($\Phi,\bD_d$)=$\pi/4$, corresponding to Case $\mathcal{C}$.\ref{case_mixed}, and size of change as given on the x-axis (multiplied by $\sqrt{d}$).}\label{f:degeneracy}
\end{center}
\end{figure}

	The following interpretation also explains the above mathematical findings: In situation { $\mathcal{C}$.\ref{case_mixed}}, each component has a common factor $\{\eta_t\}$ weighted according to $\bPhi_d$ plus some independent noise. If a change occurs in sync with the common factor it will be difficult to detect as in order to get the correct size, we have to allow for the random movements of $\{\eta_t\}$ thus increasing the critical values in that direction. In directions orthogonal to it, we only have to take the independent noise into account which yields comparably smaller noise in the projection.
	In an economic setting, this driving factor could for example be thought of as an economic factor behind certain companies (e.g.\ { ones in the same industry}). If a change occurs in those companies proportional to this driving factor it will be difficult to distinguish a different economic state of this driving factor from a mean change that is proportional to the influence of this factor.
	
\end{case}

\section{ High dimensional efficiency of panel change point test}\label{section_comparison}
In this section, we will compare the power of the above projection tests with corresponding CUSUM tests that take the full multivariate information into account.
First statistics of this type were developed for the multivariate setting with $d$ fixed \cite{horetal99}. The multivariate change point statistic (using the full multivariate information and no additional knowledge about the change) for the at most one mean change is given as a weighted maximum or sum of the following quadratic form
\begin{align}
	V^M_{d}(x)=\bZ_{T}(x)^T \boldsymbol A \bZ_{T}(x)^T,
\end{align}
where $\bZ_T(x)=(Z_{T,1}(x),\ldots,Z_{T,d}(x))^T$ is defined as in
\eqref{eq_def_Z}. The usual choice is $\boldsymbol A=\Sigma^{-1}$,
where $\Sigma$ is the covariance matrix of the multivariate observations. The weighting with $\Sigma^{-1}$ has the advantages that it (a) leads to a pivotal limit and (b) the statistic can detect all changes no matter what the direction. The second remains true for any positive definite matrix $\mathbf{A}$, the first also remains true for lower rank matrices with a decorrelation property of the errors, where this latter approach is essentially a projection (into a lower-dimensional space) as discussed in the previous sections. For an extensive discussion of this issue for the example of changes in the autoregressive structure of time series we refer to \citet{KirchMO2014}. The choice $A=\Sigma^{-1}$ corresponds to the correctly scaled case, while the misscaled case corresponds to the choice $\boldsymbol{A}=\boldsymbol M^{-1}$.

However, this multivariate setup is not very suitable for  the theoretic power comparison we are interested in because the limit distribution (a sum of $d$ squared Brownian bridges with covariance matrix $\Sigma^{1/2}\boldsymbol A\Sigma^{1/2}$) still depends on $d$ as well as the possible misspecification. Therefore, a comparison needs to take both the rates, the additive term and the noise level (which depends also on the misspecification of the covariance) present in the limit distribution into account.
For the panel data settings on the other hand, where $d\to\infty$, all the information about $d$ is  contained only in the rates rather than the limit distribution as in the previous sections. This makes the results interpretable in terms of the high dimensional efficiency. The panel null limit distribution differs from the one obtained for the projections but they are at least on the same scale, and not dependent on $d$ nor the covariance structure $\Sigma$.
Furthermore, the panel statistic is strongly related to the multivariate statistic so that the same qualitative statements can be expected, which is confirmed by simulations (results not shown).

We will now introduce the statistic for detecting changes in the mean introduced by \citet{HorvathH2012}, termed the HH Panel Statistic in this paper. Unlike in the above theory for projections, it is necessary to assume independence between components. Because the proofs are based on a central limit theorem across components, they cannot be generalized to uncorrelated (but dependent) data. For this reason, we cannot easily derive the asymptotic theory after standardization of the data. This is different from the multivariate situation, where this can easily be achieved.

We are interested in a comparison of the high dimensional efficiency for correctly specified covariance, i.e.\ $\boldsymbol A=\Sigma^{-1}$,  in addition to a comparison in the misspecified case, $\boldsymbol A=\boldsymbol M^{-1}$. The latter has  already been discussed by \citet{HorvathH2012} to some extent.  To be precise, a common factor  is introduced as in $\mathcal{C}$.3 and  the limit of the statistic (with $\boldsymbol A=\Lambda^{-1}$) under the assumption that the components are independent (i.e.\ $\Lambda$ being a diagonal matrix) is considered.
Because of the necessity to estimate the unknown covariance structure for practical purposes, the same qualitative effects as discussed here can be expected if a statistic and corresponding limit distribution were available for the covariance matrix $\Sigma$. 

\subsection{ Efficiency for panel change point tests for independent panels}
{ The above multivariate statistics have been adapted to the panel data setup under the assumption of independent components by \citet{Bai2010} for estimation as well as \citet{HorvathH2012} for testing}. Those statistics are obtained as weighted maxima or sum of the following (univariate) partial sum process
\begin{align}\label{eq_stat_panel}
	&{V}_{d,T}(x)=\frac{1}{\sqrt{d}}\sum_{i=1}^d\left( \frac{1}{\sigma_i^2}Z_{T,i}^2(x)-\frac{\lfloor Tx\rfloor (T-\lfloor Tx\rfloor)}{T^2} \right),
\end{align}
where $Z_{T,i}$ is as in \eqref{eq_def_Z} and $\sigma_i^2={ \var e_{i,1}}$.

The following theorem gives a central limit theorem for this partial sum process (under the null) from which null asymptotics of the corresponding statistics can be derived. It was proven by ~\citet[Theorem 1]{HorvathH2012}, under somewhat more general assumptions allowing in particular for time series errors (in the form of linear processes). While this makes estimation of the covariances more difficult and less precise as long-run covariances need to be estimated, it has no effect on the high dimensional efficiency. Therefore, we will concentrate on the i.i.d. (across time) situation in this work to keep things simpler { purely in terms of the calculations}.

\begin{theorem}\label{th_panel_null}
	Let Model \eqref{eq_model} hold with $\{e_{i,t}:i,t\}$ independent (where the important assumption is the independence across components) such that $\var e_{i,t}\gs c >0$ for all $i$ and $\lim\sup_{d\to\infty}\frac 1 d\sum_{i=1}^d\E |e_{i,t}|^{\nu}<\infty$ for some $\nu>4$. Furthermore, let $\frac{d}{T^2}\to 0$.  Then, it holds under the null hypothesis of no change
\begin{align*}
	&{V}_{d,T}(x)\overset{D[0,1]}{\longrightarrow}\sqrt{2} (1-x)^2W\left( \frac{x^2}{(1-x)^2} \right),
\end{align*}
where $W(\cdot)$ is a standard Wiener process.
\end{theorem}

The following theorem derives the high dimensional efficiency in this setting for HH Panel statistics such as $\max_{0\ls x\ls 1}V_{d,T}(x)$, which we use in simulations with both known and estimated standard deviations or $\int_0^1V_{d,T}(x)$.
\begin{theorem}\label{th_ind_cont_2}
Let the assumptions of Theorem~\ref{th_panel_null} on the errors be fulfilled, which implies in particular that $\Sigma=\mbox{diag}(\sigma_1^2,\ldots,\sigma_d^2)$. Then, the high dimensional efficiency of HH Panel statistic tests is given by
\begin{align*}
	\HDE_2(\bD_d)=	\frac{1}{d^{1/4}} \,\|\Sigma^{-1/2}\bD_d\|.
\end{align*}
\end{theorem}

Equivalent assertions to Corollary~\ref{th_power_1_AMOC} can be obtained analogously.

Comparing  this high dimensional efficiency with the ones given in Theorem~\ref{th_ind_cont}, Proposition~\ref{lemma_oracle_representation} as well as Theorem~\ref{th_random_proj_scaled}, we note that the high dimensional efficiency of the above HH Panel Statistic is an order ${d}^{1/4}$ worse than for the oracle but a ${d}^{1/4}$ better than the scaled random projection (also see Figure \ref{f:d_increasing}). By Theorem~\ref{th_ind_cont} we also get an indication of how wrong our assumption on $\bD_d$ can be while still obtaining a better efficiency than with the full multivariate information. More precisely, the theorems define a cone around the search direction such that the projection statistic has higher efficiency than the panel statistic if the true change direction is in this cone. 
 We can see the finite sample nature of this phenomena in Figure \ref{f:angle}.

%

 \subsection{ Efficiency of panel change point tests under dependence between Components}\label{section_panel_dep}
We now turn again to the misspecified situation, where we use the above statistic in a situation where components are not uncorrelated. Following \citet{HorvathH2012}, we consider  the mixed case $\mathcal{C}$.3 for illustration.
			The next proposition derives the null limit distribution for that special case. It turns out that the limit as well as convergence rates depend on the strength of the contamination by the common factor.
			\begin{lemma}\label{th_panel_null_miss}
				Let Case~$\mathcal{C}$.\ref{case_mixed} hold with $\nu>4$, $0<c\ls s_i\ls C<\infty$ and $\Phi_i^2\ls C<\infty$ for all $i$ and some constants $c,C$  and consider $V_{d,T}(x)$  defined as in \eqref{eq_stat_panel}, where $\sigma_i^2={ \var e_{i,1}}$ but the rest of the dependency structure is not taken into account. The asymptotic behavior of $V_{d,T}(x)$ then depends on the behavior of
				\begin{align*}
					A_d:=\sum_{i=1}^d\frac{\Phi_i^2}{\sigma_i^2}.
				\end{align*}
				\begin{enumerate}[a)]
					\item If $A_d/\sqrt{d}\to 0$, then the dependency is negligible, i.e.\
\begin{align*}
	&{V}_{d,T}(x)\overset{D[0,1]}{\longrightarrow}\sqrt{2} (1-x)^2W\left( \frac{x^2}{(1-x)^2} \right),
\end{align*}
where $W(\cdot)$ is a standard Wiener process.
					\item If $A_d/\sqrt{d}\to \xi$, $0<\xi<1$, then
\begin{align*}
	&{V}_{d,T}(x)\overset{D[0,1]}{\longrightarrow}\sqrt{2} (1-x)^2W\left( \frac{x^2}{(1-x)^2} \right)+\xi \,(B^2(x)-x(1-x)),
\end{align*}
where $W(\cdot)$ is a standard Wiener process and $B(\cdot)$ is a standard Brownian bridge.
					\item If $A_d/\sqrt{d}\to \infty$, then
						\begin{align*}
							\frac{\sqrt{d}\,V_{d,T}(x)}{A_d}\overset{D[0,1]}{\longrightarrow} B^2(x)-x(1-x),
						\end{align*}
						where $\{B(x):0\ls x\ls 1\}$ is a standard Brownian bridge.
				\end{enumerate}
\end{lemma}

Because $A_d$ in the above theorem cannot feasibly be estimated, this result cannot be used to derive critical values for panel test statistics.
Consequently, the exact shape of the limit distribution in the above lemma is not important. However, the lemma is necessary to derive the high dimensional efficiency of the panel statistics in this misspecified case.
Furthermore, it indicates that using the limit distribution from the previous section to derive critical values will result in asymptotically wrong sizes if a stronger contamination by a common factor is present.
The simulations in Figure \ref{f:size} also confirm this fact and show that the size distortion can be enormous. It does not matter whether the variance of the components in the panel statistic takes into account the dependency or simply uses the noise variance (Figure \ref{f:size}(a)), or whether a change is accounted for or not in the estimation (Figure \ref{f:size}(b)-(c)). This illustrates, that the full panel statistic is very sensitive with respect to deviations from the assumed underlying covariance structure in terms of size.

In the situation of a) and b) above, the dependency structure introduced by the common factor is still small enough asymptotically to not change the high dimensional efficiency as given in Theorem~\ref{th_ind_cont_2}, which is analogous to the proof of Theorem~\ref{th_ind_cont_2}. Therefore, we will now concentrate on situation c) in the below proposition, which is the case where the noise coming from the common factor does not disappear asymptotically.
\begin{theorem}\label{th_miss_cont}
	Let the assumptions of Lemma~\ref{th_panel_null_miss} on the errors be fulfilled and $A_d/\sqrt{d}\rightarrow\infty$, then the corresponding panel statistics have high dimensional efficiency
\begin{align*}
	\HDE_3(\bD_d)=	\frac{1}{\sqrt{A_d}} \,\sqrt{\bD_d^T\;\mbox{diag}\left(\frac{1}{s_1^2+\Phi_1^2},\ldots,\frac{1}{s_d^2+\Phi_d^2}\right)\,\bD_d}.
\end{align*}
%
\end{theorem}

The next corollary shows that the efficiency of the quasi oracle (which is scaled with $\mbox{diag}\left(\frac{1}{s_1^2+\Phi_1^2},\ldots,\frac{1}{s_d^2+\Phi_d^2}\right)$ analogously to the panel statistic) is always at least as good  as the efficiency of the panel statistic. Additionally, the efficiency of the panel statistic becomes as bad as the efficiency of the corresponding (diagonally) scaled random projection if $A_d/d\to A>0$, which is typically the case if the dependency is non-sparse and non-negligible.
\begin{cor}\label{cor_miss_panel}
	Let the assumptions of Lemma \ref{th_panel_null_miss} on the errors be fulfilled, then the following assertions hold:
	\begin{enumerate}[a)]
		\item The high dimensional efficiency of the quasi-oracle is always at least as good as the one of the misspecified panel statistic, i.e. with $\Sigma=\mbox{diag}(\sigma_1^2,\ldots,\sigma_d^2)+\bPhi\bPhi^T$, $\Lambda_d=\mbox{diag}(\sigma_1^2,\ldots,\sigma_d^2)$, it holds
			\begin{align*}
				\HDE_1^2(\bD_d,\qo)\gs \frac{\bD_d^T\Lambda_d^{-1}\bD_d}{1+A_d},
			\end{align*}
			where equality holds iff $\bD_d\sim\bPhi$.
		\item If $A_d/d\to A>0$, then the high dimensional efficiency of the panel statistic is as bad as a randomly scaled projection, i.e.
			\begin{align*}
				\HDE_3^2(\bD_d)=\frac{ \bD_d^T\Lambda_d^{-1}\bD_d}{d} \,(A_d+o(1)).
			\end{align*}
	\end{enumerate}
	\end{cor}

	In particular, for $A_d/d\to A>0$ the efficiency of the misscaled panel statistic is always as bad as the efficiency of the random projection, this only holds for the misscaled (quasi-)  projection if $\bD_d\sim \bPhi$. This effect can be clearly see in Figures \ref{f:HHcomp2} and \ref{f:degeneracy}, where in all cases H\&H Sigma refers to the panel statistic using known variance, and H\&H Var uses an estimated variance, showing again that this concept of efficiency is very well suited to understand the small sample power behavior of the corresponding statistics.


\section{Data Example}\label{sec_data}

As an illustrative example which shows the small sample behaviour of the statistics illustrated above also apply in real data, we examine the stability of change points detected by different methods in several world stock market indices. More specifically, the Fuller Log Squared Returns \cite[p 496]{fuller1996introduction} of the FTSE, NASDAQ, DAX, NIKKEI, Hang Seng and CAC \footnote{We only use a small number of series to allow reliable estimates for the covariance to be used in the full multivariate statistic.} indices for the year 2015 were examined for change points. Tests based on the multivariate statistics using full covariance estimates, a multivariate statistic using only variance estimates (i.e., a diagonal covariance structure), a projection statistic in the average direction [1,1,1,1,1,1], and a projection statistic in the direction of European countries vs non-European countries [1,-1,1,-1,-1,1] (orthogonal to the average direction) were carried out. Given the considerable dependence between the different components, we would expect economies to likely rise and fall together, justifying the use of the former projection direction. However, we think it unlikely that there will be changes of the kind that when European markets goes up, non-European markets go down, and visa versa, so take this projection as an example of direction where no change is likely. It should be noted at this point that the multivariate statistic treats both of these alternatives as equally likely. As there are possible multiple changes points in this data, we examine stability by performing binary segmentation using the proposed tests, firstly on data from January to November 2015, and then subsequently adding in the data from December 2015.

\begin{figure}
\includegraphics[width=\textwidth]{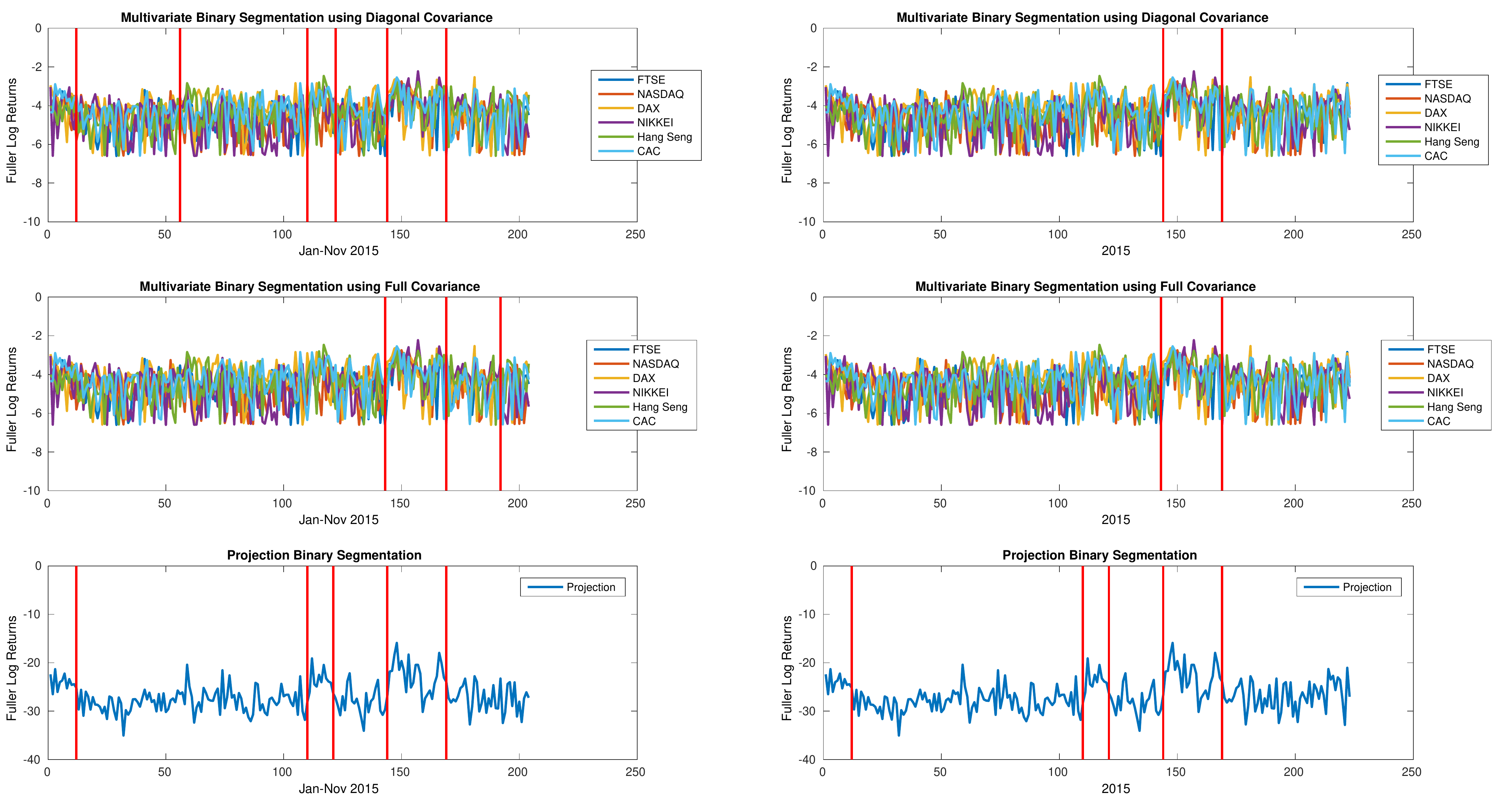}
\caption{Estimated change point locations for the market indices from binary segmentation based on different test statistics and different spans of data. First Column: Data from Jan-Nov 2015, Second Column: Data from all of 2015. First Row: Multivariate statistic with full covariance estimation; Second Row: Multivariate Statistic with Diagonal Variance Estimate; Third Row: Projection Statistic in direction [1,1,1,1,1,1]. Red vertical lines indicate changes deemed to be significant at 5\% level.}\label{f:FLSR_BinSeg}
\end{figure}

As can be seen in Figure \ref{f:FLSR_BinSeg}, the multivariate test statistic is considerably less robust than the average projection based statistic, both to the length of the data, as well as to the choice of the covariance estimate. The major cause of this instability was that the CUSUM statistic over time had two peaks, but the location of the maximal peak differed from one to the other when further data was added. This caused knock-on effects in the entire binary segmentation. Here, in all cases, independence in time was assumed as once the changes were accounted for, there was little evidence of temporal dependence in the data. However, even if time series dependence is accounted for by using an estimate of the long run covariance in place of the independent covariance estimate, there is no difference in the qualitative conclusions (although the change points themselves varied considerably in all cases depending on the parameters chosen in the long run covariance estimation procedure \cite{politis05}). In addition, the projection estimate was robust to whether the direction was scaled by the full covariance, the diagonal of the covariance or not scaled at all, as well as to increasing the length of the data.

\begin{table}
\centering
\caption{Location, Statistic and p-value for the changes found in the 2015 market index data. (Limit Distributions: Multivariate: sum of six independent Brownian Bridges; Projection: Single Brownian Bridge)}\label{t:FLSR_pvals}
\begin{tabular}{llllll}
\hline
\multicolumn{6}{c}{Multivariate: Full Covariance}\\
\hline
Day&&&&143&169\\
Statistic Value&&&&6.8541&5.1581\\
p&&&&0.0012&0.0173\\
\hline
\multicolumn{6}{c}{Multivariate: Diagonal Covariance}\\
\hline
Day&&&&144&169\\
Statistic Val&&&&9.9995&11.7030\\
p-value&&&&$<10^{-5}$&$<10^{-5}$\\
\hline
\multicolumn{6}{c}{Projection: scaled [1 1 1 1 1 1]}\\
\hline
Day&12&110&121&144&169\\
Statistic Value&2.1307&3.5390&2.9518&3.3173&2.0900\\
p&0.0285&0.0017&0.0057&0.0027&0.0307\\
\hline
\end{tabular}
\end{table}

The p-values for the changes on full year's data are given in Table \ref{t:FLSR_pvals}. While it can be seen that the projection p-vals are larger for the two common change points than in the multivariate case, the same changes are detected with all methods. However, additional changes are found with the projection method, and the p-vals are well below the critical value of 5\%. This shows that having knowledge of the likely direction of change can allow further changes to be found beyond those in an unrestricted multivariate search. As expected though, using an unlikely direction does not find change points, with the hypothesis that there are no changes which affect European markets differently to non-European markets being accepted (p=0.18).

\section{Conclusions}\label{section_conclusions}

	The primary aims of this paper were to introduce a theoretic method to compare the small sample behavior of different high dimensional tests by asymptotic methods. Furthermore, projection based statistics are introduced into the analysis of change points in high dimensions and compared and contrasted with the panel based statistics that are currently available.
	The new concept of high dimensional efficiency allows a comparison of the magnitude of changes that can be detected asymptotically as the number of dimensions increases. All the tests in the paper are benchmarked against random projections. Because the space covered by far away angles increases rapidly with the dimension, the power of these becomes very poor in higher dimensions rendering random projections useless in practice { for detecting change points}.

	In summary, the following two assertions were proven:
First, a suitable projection will substantially increase the power of detection but at the cost of a loss in power if the change is at a large angle away from the projection vector. Second, projections are more robust compared to the panel based statistic with respect to misspecification in the covariance structure both in terms of size and power.

The panel statistic \cite{Bai2010,HorvathH2012} test works well in situations where the panels are independent across dimension, and there is little to no information about the direction of the change. However, as soon as dependency is present, the size properties of these statistics become difficult and their high dimensional efficiencies mimic those of random projections. Misspecification of the covariance structure can be problematic for all tests. However, if the direction of the likely change is known, then it is always preferable to use the corresponding projection (scaled with the assumed covariance structure), rather than either the panel statistic or a random projection, regardless of whether the covariance is misspecified or not.

This results in this paper raise many questions for future work. It would be of considerable interest to determine whether projections can be derived using data driven techniques, such as sparse PCA, for example, and whether such projections would be better than random projections. Preliminary work suggests that this may be so in some situations but not others, and a nice procedure by \citet{WangS2016} investigates a related technique. Further many multiple change point procedures use binary segmentation or related methods to find the multiple change points, so much of the work here would apply equally in suitably defined sub intervals which are then assumed to contain at most one change. This was the approach taken in the data example here. In addition, all the results here have been assessed with respect to choosing a single projection for the test which is optimal if the direction of the change is known. However, in some situations only qualitative information is known or several change scenarios are of interest. Then, it could be very beneficial to determine how best to combine this information into testing procedures based on several projections, where a standard subspace approach may not be ideal as the information about the likely direction of changes is lost. 
Finally, while the framework in this paper concentrates on tests with a given size, as soon as a-priori information is considered, then it is natural to ask whether related Bayesian approaches are of use, and indeed quantifying not only the a-priori direction of change, but also its uncertainty, prior to conducting the test is a natural line of further research.

\section{Proofs}\label{section_proofs}
\begin{proof}[of Theorem~\ref{th_proj_null}]
	We need to prove the following functional central limit theorem for the triangular array of projected random variables $Y_{t,d}=\sum_{j=1}^dp_j(d)e_{j,t}(d)$ given the (possibly random) projection $\bp_d=(p_1(d),\ldots,p_d(d))^T$:
	\begin{align}\label{eq_fclt_proj}
		&\left\{\frac{1}{\sqrt{T \tau^2(\bp_d)}}\sum_{t=1}^{\lfloor Tx\rfloor}Y_{t,d}:0\ls x\ls 1\,|\,\bp_d\right\}\overset{D[0,1]}{\longrightarrow} \{W(x):0\ls x\ls 1\}\qquad a.s.,
	\end{align}
	where $\{W(\cdot)\}$ denotes a standard Wiener process.
	
	The proof for tightness is analogous to the one given in Theorem 16.1 of \citet{Billingsley1968} as it only depends on the independence across time (which also holds conditionally given $\bp_d$ due to the independence of $\bp_d$ and $\{\be_t(d)\}$). Similarly, the proof for the convergence of the finite dimensional distributions follows the proof of Theorem 10.1 in \citet{Billingsley1968}, where we need to use the Lindeberg-Levy-version of the univariate central limit theorem for triangular arrays. More precisely, we need to prove the  Lindeberg condition given by
	\begin{align*}
		&\E \left( \frac{Y_{1,d}^2}{\tau^2(\bp_d)}1_{\{Y_{1,d}/\tau(\bp_d)\gs \eps \sqrt{T}\}} \,|\,\bp_d\right)\to 0 \quad a.s.
	\end{align*}
	for any $\eps>0$. The following Lyapunov-type condition implies the above Lindeberg condition:
	\begin{align}\label{eq_lyap}
		\E\left(\left|\frac{Y_{1,d}}{\tau(\bp_d)} \right|^{\nu}\,|\,\bp_d\right)=\E\left(\left|\frac{\bp_d^T\be_1(d)}{\tau(\bp_d)}\right|^{\nu}\,|\,\bp_d\right)=o(T^{\nu/2-1})\qquad a.s.,
	\end{align}
	where  $\nu>2$ as given in the theorem. 
	Let
	\begin{align*}
		\tilde{\bp}_d=\frac{\bp_d}{\sqrt{\bp_d^T\cov \be_1(d)\bp_d}},
	\end{align*}
	then the above Lyapunov condition is equal to
	\begin{align*}
		&\E\left( \left|\tilde{\bp}_d^T\be_1(d)\right|^{\nu}\,|\,\bp_d\right)=o(T^{\nu/2-1})\qquad a.s.
	\end{align*}
	In the situation of a) $\cov\be_1(d)=\sum_{j\gs 1}\ba_j(d)\ba_j^T(d)$ and we get by the Rosenthal inequality (confer e.g.\  \citet[9.7c]{LinB2010})
\begin{align*}
	&	\E\left(\left|\sum_{j=m}^n\tilde{\bp}_d^T\ba_{j}(d)\eta_{j,1}(d) \right|^{\nu}\,|\,\bp_d\right)\\
	&\ls  O(1)\sum_{j=m}^n\left|\tilde{\bp}_d^T\ba_j(d)\right|^{\nu}\;\E|\eta_{j,1}(d)|^{\nu}+ O(1) \left( \sum_{j=m}^n\left(\tilde{\bp}_d^T\ba_j(d)\right)^2\,\var\eta_{j,1}(d) \right)^{\nu/2},\end{align*}
where   the right-hand side is bounded for any $m,n$ with a bound that does not depend on $T$ or $d$ and converges to zero for $m,n\to\infty$ as $\E |\eta_j(d)|^{\nu}\ls C$ hence $\var\eta_j(d)\ls 1+C$ and by definition of $\tilde{\bp}_d$ it holds $\sum_{j=m}^n|\tilde{\bp}_d^T\ba_j(d)|^2\ls \tilde{\bp}_d^T\cov \be_1(d)\,\tilde{\bp}_d\ls 1$, hence also
$|\tilde{\bp}_d^T\ba_j(d)|^{\nu}\ls |\tilde{\bp}_d^T\ba_j(d)|^2$ and  $\sum_{j=m}^n|\tilde{\bp}_d^T\ba_j(d)|^{\nu}\ls 1$.

Consequently, the infinite series exists in an $L^{\nu}$-sense with the following uniform (in $T$ and $d$) moment  bound
\begin{align}\label{eq_lyap1}
	\E \left(\left|\tilde{\bp}_d^T\be_1(d)\right|^{\nu}\,|\,\bp_d\right)=O(1)=o(T^{\nu/2-1})\qquad a.s.
\end{align}

To prove the Lyapunov-condition under the assumptions of b) we use the Jenssen-inequality which yields
\begin{align}\label{eq_lyap2}
	&\E\left( \left|\tilde{\bp}_d^T\be_1(d)\right|^{\nu}\,|\,\bp_d\right) = \|\tilde{\bp}_d\|_1^{\nu}\,\E \left(\left( \sum_{i=1}^d\frac{|\tilde{p}_{i,d}|}{\|\tilde{\bp}_d\|_1}\,|e_{i,1}(d)| \right)^{\nu}\,|\,\bp_d\right)\notag\\
	&
	\ls  \|\tilde{\bp}_d\|_1^{\nu} \sum_{i=1}^d\frac{|\tilde{p}_{i,d}|}{\|\tilde{\bp}_d\|_1}\,\E|e_{i,1}(d)|^{\nu}\ls C 	 \left(\frac{\|\bp_d\|_1}{\sqrt{\bp_d^{T}\cov (\be_1(d))\bp_d^T}}\right)^{\nu}=o(T^{\nu/2-1})\qquad a.s
.
\end{align}

\end{proof}

\begin{proof}[of Lemma~\ref{lem_variance}]
With the notation of the proof of Theorem~\ref{th_proj_null}
both estimators (as functions of $\bp_d$) fulfill ($j=1,2$)
	\begin{align*}
		\frac{\widehat{\tau}_{j,d,T}^2(\bp_d)}{\tau^2(\bp_d)}=\widehat{\tau}_{j,d,T}^2(\tilde{\bp}_d).
	\end{align*}
	First by the independence across time we get by the van Bahr-Esseen inequality (confer e.g. \citet[9.3 and 9.4]{LinB2010}) for some constant $C>0$, which may differ from line to line,
	\begin{align}\label{eq_sum_proof_test}
		&	\E_{\bp_d}\left|\sum_{j=a+1}^b\left(\left(\tilde{\bp}_d^T\be_j(d)\right)^2-1\right)\right|^{\nu/2}\ls C (b-a)^{\max(1,\nu/4)}\,\E_{\bp_d}\left|\left(\tilde{\bp}_d^T\be_1(d)\right)^2-1\right|^{\nu/2}\notag\\
		&\ls C\,(b-a)^{\max(1,\nu/4)}\,\max\left(1, \E_{\bp_d}\left|\tilde{\bp}_d^T\be_1(d)\right|^{\nu}\right)\notag\\
	&\ls \begin{cases}
		C (b-a)^{\max(1,\nu/4)}\qquad a.s., & \text{in  a)},\\
		C (b-a)^{\max(1,\nu/4)} \max\left( 1,\left( \frac{\|\bp_d\|_1}{\sqrt{\bp_d^T\cov \be_1(d)\bp_d}} \right)^{\nu} \right), &\text{in  b)},
	\end{cases}
\end{align}
by \eqref{eq_lyap1} resp.\ \eqref{eq_lyap2},
where $\E_{\bp_d}$ denotes the conditional expectation given $\bp_d$.
An application of the Markov-inequality now yields for any $\epsilon>0$
\begin{align*}
	&P\left(\frac 1 T\left| \sum_{j=1}^T\left(\left(\tilde{\bp}_d^T\be_j(d)\right)^2-1\right) \right|\gs \epsilon \,\Big |\, \bp_d \right)\\
	&\ls
	\begin{cases} \frac{C}{\epsilon^{\nu/2}}T^{-\nu/2+\max(1,\nu/4)}\qquad a.s., &\text{in a)},\\
		\frac{C}{\epsilon^{\nu/2}}T^{-\nu/2+\max(1,\nu/4)}o(T^{\nu/2-\nu/\min(\nu,4)})\qquad a.s., &\text{in b)},
	\end{cases}\\
	&
	\to 0 \qquad a.s.
\end{align*}
Similar arguments yield
\begin{align*}
	&P\left(\frac 1 T\left|\sum_{j=1}^T\tilde{\bp}_d^T\be_j(d)\right|\gs \eps \,\Big|\, \bp_d\right)\to 0\qquad a.s.
\end{align*}
proving a) and b) for $\widehat{\tau}_{1,d,T}^2(\bp_d)$.

From \eqref{eq_sum_proof_test} it follows by Theorem B.1 resp. B.4 in \citet{kirchdiss}
\begin{align*}
	&	\E_{\bp_d}\max_{1\ls k\ls T}\left| \sum_{j=1}^k \left(\left(\tilde{\bp}_d^T\be_j(d)\right)^2-1\right)\right|^{\nu/2}\\
	&\ls \begin{cases}
		C T^{\max(1,\nu/4)}(\log T)^{\frac{(4-\nu)_+\nu}{2(4-\nu)}}\qquad a.s., & \text{in a)},\notag\\
		C T^{\max(1,\nu/4)} (\log T)^{\frac{(4-\nu)_+\nu}{2(4-\nu)}}\max\left( 1,\left( \frac{\|\bp_d\|_1}{\sqrt{\bp_d^T\cov \be_1(d)\bp_d}} \right)^{\nu} \right), &\text{in b)},
	\end{cases}
	\\
	&\to 0\qquad a.s.
\end{align*}
An application of the Markov inequality now yields for any $\eps>0$
\begin{align*}
	&P\left(\max_{1\ls k\ls T}\frac{1}{T} \left|\sum_{j=1}^k\left((\left(\tilde{\bp}_d^T\be_j(d)\right)^2-1\right)\right|\gs \eps\,\Big|\, \bp_d\right)\to 0\qquad a.s.\end{align*}
By the independence across time it holds
\begin{align*}
	\left\{\sum_{j=k+1}^T\left(\left(\tilde{\bp}_d^T\be_j(d)\right)^2-1\right):1\ls k\ls T\right\}\deq 	 \left\{\sum_{j=1}^{T-k}\left(\left(\tilde{\bp}_d^T\be_j(d)\right)^2-1\right):1\ls k\ls T\right\},
\end{align*}
which implies
\begin{align*}
	&P\left(\max_{1\ls k\ls T}\frac{1}{T} \left|\sum_{j=k+1}^T\left(\left(\tilde{\bp}_d^T\be_j(d)\right)^2-1\right)\right|\gs \eps\,\Big|\, \bp_d\right)\to 0\qquad a.s.
\end{align*}
Similar assertions can be obtained along the same lines for $\max_{1\ls k\ls T}\frac{1}{T} \left|\sum_{j=1}^k\tilde{\bp}_d^T\be_j(d)\right|$ as well as $\max_{1\ls k\ls T}\frac{1}{T}\left| \sum_{j=k+1}^T\tilde{\bp}_d^T\be_j(d)\right|$, which imply the assertion for $\widehat{\tau}_{2,d,T}^2(\bp_d)$.
\end{proof}

\begin{proof}[of Corollary~\ref{th_asym_null_stat}]
	By an  application of the continuous mapping theorem and Theorem~\ref{th_proj_null} we get the assertions for the truncated maxima resp.\ the sums over $[\tau T, (1-\tau) T]$ for any $\tau>0$ towards equivalently truncated limit distributions.
	Because we assume independence across time (with existing second moments) the H\'ajek-R\'enyi inequality yields for all $\eps>0$	
	\begin{align*}
	&P\left(	\max_{1\ls k\ls \tau T}w(k/T)\left|\sum_{t=1}^k\tilde{\bp}^T_d\be_t(d)\right|\gs \eps\,\Big|\,\bp_d\right)\to 0\quad a.s.\\
	&P\left(\max_{(1-\tau) T\ls k\ls }w(k/T)\left|\sum_{t=k+1}^T\tilde{\bp}^T_d\be_t(d)\right|\gs \eps\,\Big|\,\bp_d\right)\to 0\quad a.s.
	\end{align*}
	as $\tau\to 0$ uniformly in $T$, where the notation of the proof of Theorem~\ref{th_proj_null} has been used.
	This in addition to an equivalent argument for the limit process shows that the truncation is asymptotically negligible proving the desired results.
\end{proof}

\begin{proof}[of Theorem~\ref{th_ind_cont}]
		We consider the situation where $\sqrt{T}\,\HDE_1(\bD_d,\bp_d)$ converges a.s.
	Under alternatives it holds
\begin{align*}
	\frac{U_{d,T}(x)}{\tau(\bp_d)}=\frac{U_{d,T}(x;\be)}{\tau(\bp_d)}+\mbox{sgn}(\bD_d^T\bp_d)\,\sqrt{T}\, \HDE_1(\bD_d,\bp_d)\left( \frac 1 T \sum_{i=1}^{\lfloor T x\rfloor}g(i/T)-\frac{\lfloor T x\rfloor}{T^2}\sum_{j=1}^Tg(j/T) \right),
\end{align*}
where $U_{d,T}(x;\be)$ is the corresponding functional of the error process. By Theorem~\ref{th_proj_null} it holds
\begin{align*}
	\left\{\frac{U_{d,T}(x;\be)}{\tau(\bp_d)}:0\ls x\ls 1\,|\,\bp_d\right\}\overset{D[0,1]}{\longrightarrow} \{B(x):0\ls x\ls 1\}\qquad a.s.
\end{align*}
Furthermore, by the Riemann-integrability of $g(\cdot)$ it follows
\begin{align*}
	\sup_{0\ls x\ls 1}\left|	 \frac 1 T \sum_{i=1}^{\lfloor T x\rfloor}g(i/T)-\frac{\lfloor T x\rfloor}{T^2}\sum_{j=1}^Tg(j/T) -\left(\int_0^xg(t)\,dt-x\int_0^1g(t)\,dt\right)\right|\to 0.
\end{align*}
For any $\tau >0$
	\begin{align*}
		&\max_{ \tau \ls k/T \ls 1-\tau} w^2(k/T) \frac{U_{d,T}^2(k/T)}{\tau^2(\bp_d)}\\
		&=T\,\HDE^2_1(\bD_d,\bp_d)\left(\sup_{\tau\ls x\ls 1-\tau}w^2(x)\,\left( \int_0^xg(t)\,dt-x\int_0^1g(t)\,dt \right)^2+o_{P_{\bp_d}}(1)\right) \qquad a.s.,
	\end{align*}
	where $P_{\bp_d}$ denotes the conditional probability given $\bp_d$. Because by assumption $\sup_{\tau\ls x\ls 1-\tau}w^2(x)\,\left( \int_0^xg(t)\,dt-x\int_0^1g(t)\,dt \right)^2>0$ for some $\tau>0$, so that the above term becomes unbounded asymptotically. This gives the assertion for the max statistics, similar arguments give the assertion for the sum statistic.
\end{proof}

\begin{proof}[of Corollary~\ref{th_power_1_AMOC}]
	Similarly to the proof of Theorem~\ref{th_ind_cont} it follows (where the uniformity at $0$ and $1$ follows by the assumptions on the rate of divergence for $w(\cdot)$ at $0$ or $1$)
	\begin{align*}
		\sup_{0< x< 1}w^2(x)\,\left|\frac{U^2_{d,T}(x)}{\tau^2(\bp_d)T\,\HDE^2_1(\bD_d,\bp_d)}-\left( (x-\vth)_+-x(1-\vth) \right)^2\right|=o_{P_{\bp_d}}(1)\qquad a.s.,
	\end{align*}
	which implies the  assertion by standard arguments on noting that
	\begin{align*}
		\widehat{\vth}_T=\arg\max_{0\ls x\ls 1} w^2(x)\,\frac{U^2_{d,T}(x)}{\tau^2(\bp_d)T\,\HDE^2_1(\bD_d,\bp_d)},\quad \vth=\arg\max_{0\ls x\ls 1} w^2(x)\,\left( (x-\vth)_+-x(1-\vth) \right)^2.
	\end{align*}
\end{proof}

\begin{proof}[of Proposition~\ref{lemma_oracle_representation}]
The assertion follows from
\begin{align*}
	&\tau^2(\bp_d)=\bp_d^T\Sigma\bp_d=\|\Sigma^{1/2}\bp_d\|^2,\\
	&|\langle \bD_d,\bp_d\rangle|=(\Sigma^{-1/2}\bD_d)^T(\Sigma^{1/2}\bp_d)=\|\Sigma^{-1/2}\bD_d\|\,\|\Sigma^{1/2}\bp_d\|\,\cos(\alpha_{\Sigma^{-1/2}\bD_d,\Sigma^{1/2}\bp_d}).
\end{align*}
\end{proof}

\begin{proof}[of Theorem~\ref{th_random_proj_scaled}]
	Let $\boldsymbol X_d=(X_1,\ldots,X_d)^T$ be N(0,$I_d$), then by \citet{Marsaglia72}  it holds $\br_d\deq (X_1,\ldots,X_d)^T/\|(X_1,\ldots,X_d)^T\|$ and it follows by \eqref{eq_22}
\begin{align*}
	& \HDE_1^2(\bD_d,\Sigma^{-1/2}\br_d)\,\frac{d}{\|\Sigma^{-1/2}\bD_d\|^2}
	\deq\frac{\left|\frac{\boldsymbol X_d^T\Sigma^{-1/2}\,\bD_d}{\|\Sigma^{-1/2}\bD_d\|}\right|^2}{\frac{\boldsymbol X_d^T\boldsymbol X_d}{\E\boldsymbol X_d^T\boldsymbol X_d }}
\end{align*}
Since the numerator has a $\chi^2_1$ distribution (not depending on $d$), there exist for any $\eps>0$  constants $0<c_1<C_1<\infty$ such that
\begin{align*}
	&\sup_{d\gs 1} P\left(c_1\ls \left|\frac{\boldsymbol X_d^T\Sigma^{-1/2}\,\bD_d}{\|\Sigma^{-1/2}\bD_d\|}\right|^2\ls C_1\right)\gs 1-\eps.
\end{align*}
Furthermore, the denominator has a $\chi^2_d$-distribution divided by its expectation, consequently an application of the Markov-inequality yields for any $\eps>0$ the existence of $0<C_2<\infty$ such that
\begin{align*}
	\sup_{d\gs 1}P\left( \frac{\boldsymbol X_d^T\boldsymbol X_d}{\E\boldsymbol X_d^T\boldsymbol X_d }\gs C_2
	\right)\ls \eps.
\end{align*}
By  integration by parts we get $\E \left(  \boldsymbol X_d^T\boldsymbol X_d\right)^{-1}\ls 2/d$ for $d\gs 3$ so that another application of the Markov-inequality yields that for any $\eps>0$ there exists $c_2>0$ such that
\begin{align*}
	\lim\sup_{d\to \infty}P\left( \frac{\boldsymbol X_d^T\boldsymbol X_d}{\E\boldsymbol X_d^T\boldsymbol X_d }\ls c_2
	\right)\ls \eps,
\end{align*}
completing the proof of the theorem by standard arguments.
\end{proof}

\begin{proof}[of Theorem~\ref{th_random_proj}]
	Let $\boldsymbol X_d=(X_1,\ldots,X_d)^T$ be N(0,$I_d$), then as in the proof of Theorem~\ref{th_random_proj_scaled} it holds
 \begin{align*}
	 & \HDE_1^2(\bD,\bM^{-1/2}\br_d)\,\frac{\mbox{tr}(\bM^{-1/2}\Sigma\bM^{-1/2})}{\|\bM^{-1/2}\bD_d\|^2}\deq\frac{\left|\frac{\boldsymbol X_d^T\bM^{-1/2}\,\bD_d}{\|\bM^{-1/2}\bD_d\|}\right|^2}{\frac{\boldsymbol X_d^T\bM^{-1/2}\Sigma \bM^{-1/2}\boldsymbol X_d}{\mbox{tr}(\bM^{-1/2}\Sigma\bM^{-1/2})}}.
\end{align*}
	The
	proof of the lower bound is analogous to the proof of Theorem~\ref{th_random_proj_scaled} by noting that ($A=\bM^{-1/2}\Sigma\bM^{-1/2}$)
 	\begin{align*}
	\E \bX^T A\bX=\E \sum_{i,j=1}^da_{i,j}X_iX_j=\sum_{i,j=1}^da_{i,j}\delta_{i,j}=\sum_{i=1}^da_{i,i}=\mbox{tr}(A).
\end{align*}
For the proof of the upper bound, first note that by a spectral decomposition it holds
\begin{align*}
	\frac{\bX^T\bM^{-1/2}\Sigma \bM^{-1/2}\bX}{\mbox{tr}(\bM^{-1/2}\Sigma\bM^{-1/2})}\deq \sum_{j=1}^d\alpha_jX_j^2,\qquad\text{for some } 0<\alpha_d\ls\ldots\ls \alpha_1,\quad \sum_{j=1}^d\alpha_j=1.
\end{align*}
From this we get
on the one hand by the Markov inequality
\begin{align*}
	P\left( \sum_{j=1}^d\alpha_jX_j^2\ls c \right)\ls P(\alpha_1X_1^2\ls c)\ls \left(\frac{c}{\alpha_1}\right)^{1/4}\,\E(|X_1^2|^{-1/4}),
\end{align*}
where $\E(|X_1^2|^{-1/4})=\Gamma(1/4)/(2^{1/4}\sqrt{\pi})$ exists (as can be seen using the density for a $\chi^2_1$-distribution).
On the other hand it holds for any $c\ls 1/2$ by another application of the Markov inequality
\begin{align*}
&	P\left(  \sum_{j=1}^d\alpha_jX_j^2\ls c\right)\ls P\left(\left |\sum_{j=1}^d\alpha_jX_j^2-1\right|\gs 1/2 \right)\ls 8 \sum_{i=1}^d\alpha_i^2\ls 8\alpha_1.
\end{align*}
By chosing $c=\min(1/2,(\E(|X_1^2|^{-1/4}))^{-4}/8\,\eps^5)$ we finally get
\begin{align*}
&\sup_{0<\alpha_d\ls \ldots\ls\alpha_1,\sum_{i=1}^d\alpha_i=1}	P\left( \sum_{j=1}^d\alpha_jX_j^2\ls c \right)\\
&\ls \sup_{0<\alpha_d\ls \ldots\ls\alpha_1,\sum_{i=1}^d\alpha_i=1}	\min\left(\eps \left( \frac{\eps}{8\alpha_1}\right)^{1/4} \,,\, 8\alpha_1  \right)\ls \eps,
\end{align*}
completing the proof.
\end{proof}



\begin{proof}[of Theorem~\ref{theorem_oracle_rand}]
	By the Cauchy-Schwarz inequality
	\begin{align*}
		&\tau^2(\bM^{-1}\bD_d)=\bD_d^T\bM^{-1}\sum_{j\gs 1}\ba_j\ba_j^T\bM^{-1}\bD_d=\sum_{j\gs 1}(\ba_j^T\bM^{-1}\bD_d)^2\ls \sum_{j\gs 1}\ba_j^T\bM^{-1}\ba_j\;\bD_d^T\bM^{-1}\bD_d\\
		&=\mbox{tr}\left(\bM^{-1/2}\sum_{j\gs 1}\ba_j\ba_j^T\bM^{-1/2}\right)\,\bD_d^T\bM^{-1}\bD_d, \end{align*}
	which implies  the assertion by \eqref{eq_22}.
\end{proof}

\begin{proof}[of Proposition~\ref{prop_oracle}]
Assertion a) follows from
	\begin{align*}
		&|\langle \bD_d,\po\rangle|^2=\left(\sum_{i=1}^d\frac{\delta_{i,T}^2}{\sigma_i^2}\,\sigma_i^2\right)^2\gs c^2 \left(\sum_{i=1}^d\frac{\delta_{i,T}^2}{\sigma_i^2}\right)^2=c^2 \left|\langle \bD_d,\qo\rangle\right|^2,\\
		&\tau^2(\po)=\po^T\Sigma\po=\sum_{i=1}^d\frac{\delta_{i,T}^2}{\sigma_i^2}\,\sigma_i^4\ls C^2 \left|\langle \bD_d,\qo\rangle\right|.
	\end{align*}
	Concerning b) first note that by the Cauchy-Schwarz  inequality with $\Lambda=\mbox{diag}(\sigma_1^2,\ldots,\sigma_d^2)$
\begin{align*}
	\tau^2(\qo)=\sum_{j\gs 1}(\bD_d^T\Lambda^{-1}\ba_j)^2\ls \bD_d^T\Lambda^{-2}\bD_d\,\sum_{j\gs 1}\ba_j^T\ba_j\ls \frac{\bD_d^T\bD_d}{c^2}\, \mbox{tr}(\Sigma).
\end{align*}
This implies  assertion b) by \eqref{eq_22} on noting that
\begin{align*}
	|\bD_d^T\Lambda^{-1}\bD_d|^2\gs \frac{ |\bD_d^T\bD_d|^2}{C^2}.
\end{align*}
\end{proof}

\begin{proof}[of Equation~\ref{ex_oracle}]
	By Proposition~\ref{lemma_oracle_representation} it holds for $\bD_d= k\,\bPhi_d$
\begin{align*}
	\HDE_1^2(\bD_d,\bo)=\|\Sigma^{-1/2}\bD_d\|^2= \bD_d^T(\boldsymbol D+\bPhi_d\bPhi_d^T)^{-1}\bD_d,
\end{align*}
where $\boldsymbol D=\mbox{diag}(s_1^2,\ldots,s_d^2)^T$. Hence
\begin{align*}
	&\bD_d^T(\boldsymbol D+\bPhi_d\bPhi_d^T)^{-1}\bD_d
	=(\boldsymbol D^{-1/2}\bD_d)^T\left( I_d+(\boldsymbol D^{-1/2}\bPhi_d)(\boldsymbol D^{-1/2} \bPhi_d)^T\right)^{-1}\boldsymbol{D}^{-1/2}\bD_d\\
	&=\frac{(\boldsymbol{D}^{-1/2}\bD_d)^T\boldsymbol{D}^{-1/2}\bD_d}{1+\boldsymbol{D}^{-1/2}\bPhi_d^T\boldsymbol{D}^{-1/2}\bPhi_d},
\end{align*}
where the last line follows from the fact that $\boldsymbol D^{-1/2}\bD_d=k \boldsymbol D^{-1/2}\bPhi_d$ is an eigenvector of $ I_d+(\boldsymbol D^{-1/2}\bPhi_d)(\boldsymbol D^{-1/2} \bPhi_d)^T$ with eigenvalue $1+(\boldsymbol{D}^{-1/2}\bPhi_d)^T\boldsymbol{D}^{-1/2}\bPhi_d$ hence also for the inverse of the matrix with inverse eigenvalue.
\end{proof}

\begin{proof}[of Theorem~\ref{th_ind_cont_2}]
Similarly as in the proof of Theorem~\ref{th_ind_cont} it holds
\begin{align*}
	&Z_{T,i}(x)=Z_{T,i}(x;\be)+\delta_{i,T}\sqrt{T}\left( \frac{1}{T}\sum_{j=1}^{  \lfloor Tx\rfloor}g(j/T)+\frac{\lfloor T x\rfloor}{T^2}\sum_{j=1}^Tg(j/T) \right),
\end{align*}
where $Z_{T,i}(x;\be)$ is the corresponding functional for the error sequence (rather than the actual observations). From this it follows
\begin{align*}
	&V_{d,T}(x)=V_{d,T}(x;\be)+T\,\HDE_2^2(\bD_d)\,\left( \int_0^xg(t)\,dt-x\int_0^1g(t)\,dt+o(1) \right)+R_T(x),
\end{align*}
where $R_T(x)$ is the mixed term given by
\begin{align*}
	R_T(x)=\frac{2\sqrt{T}}{\sqrt{d}}\sum_{i=1}^d\frac{\delta_{i,T}}{\sigma_i^2}Z_{T,i}(x;\be)\,\left( \int_0^xg(t)\,dt-x\int_0^1g(t)\,dt+o(1) \right)
\end{align*}
which by an application of the \haj-R\'enyi inequality (across time) yields
	\begin{align*}
		&		P\left(\sup_{0\ls x\ls 1}|R_T(x)|\gs c\right)=O\left( 1 \right) \frac{1}{c^2} T\frac 1 d \sum_{i=1}^d\frac{\delta_i^2}{\sigma_i^2}=O_P(1) \frac{1}{c^2\,\sqrt{d}}\,T\,\HDE_2(\bD_d).
	\end{align*}

From this the assertion follows by an application of Theorem~\ref{th_panel_null}.
\end{proof}

\begin{proof}[of Lemma~\ref{th_panel_null_miss}]
	The proof follows closely the proof of (28) -- (30) in \citet{HorvathH2012} but where we scale diagonally with the true variances. We will give a short sketch for the sake of completeness. The key is the following decomposition
	\begin{align*}
		&V_{d,T}(x)\\&=\frac{1}{\sqrt{d}}\sum_{i=1}^d\left( \frac{s_i^2}{s_i^2+\Phi_i^2} \frac 1 T \left( \sum_{t=1}^{\lfloor Tx\rfloor}\eta_{i,t}(d)-\frac{\lfloor T x\rfloor}{T}\sum_{t=1}^T \eta_{i,T}(d) \right)^2-\frac{\lfloor T x\rfloor \,(T-\lfloor T x\rfloor)}{T^2} \right)\\
		&+\frac{2}{\sqrt{d}}\left( \sum_{i=1}^d\frac{\Phi_i s_i}{s_i^2+\Phi_i^2}\,\frac{1}{\sqrt{T}}\left( \sum_{t=1}^{\lfloor T x\rfloor}\eta_{i,t}(d)-\frac{\lfloor Tx \rfloor}{T}\sum_{t=1}^T\eta_{i,T}(d) \right) \right)\, \frac{1}{\sqrt{T}}\left( \sum_{t=1}^{\lfloor T x\rfloor}\eta_{d+1,t}(d)-\frac{\lfloor T x\rfloor}{T}\sum_{t=1}^T\eta_{d+1,t}(d) \right)\\
		&+\frac 1 T \left( \sum_{t=1}^{\lfloor T x\rfloor}\eta_{d+1,t}(d)-\frac{\lfloor T x\rfloor}{T}\sum_{t=1}^T\eta_{d+1,t}(d) \right)^2\,\frac{1}{\sqrt{d}}A_d.
	\end{align*}
	The first term converges to the limit given in a). To see this, note that the proof of the Lyapunov condition in \citet{HorvathH2012} following equation (39) still holds because $s_i^2/(s_i^2+\Phi_i^2)$ is uniformly bounded from above by assumption (showing that the numerator is bounded) while again by assumption
\[
 \frac 1 d\sum_{i=1}^d\frac{s_i^4}{(s_i^2+\phi_i^2)^2}\gs D>0,
\]
showing that the denominator is bounded. Similarly, the proof of tightness in \citet{HorvathH2012} (equations (43) and following) remains valid. The asymptotic variance remains the same under a) and b) because by assumption
\begin{align*}
	\left| \frac 1 d \sum_{i=1}^d\frac{s_i^4}{(s_i^2+\Phi_i^2)^2}-1\right|\ls \frac{3}{d}A_d\to 0.
\end{align*}
The middle term in the above decomposition is bounded by an application of the \haj-R\'enyi inequality
\begin{align*}
&P\left(\sup_{0<x<1}\frac{1}{\sqrt{d}}\left| \sum_{i=1}^d\frac{\Phi_i s_i}{s_i^2+\Phi_i^2}\,\frac{1}{\sqrt{T}}\left( \sum_{t=1}^{\lfloor T x\rfloor}\eta_{i,t}(d)-\frac{\lfloor Tx \rfloor}{T}\sum_{t=1}^T\eta_{i,T}(d) \right) \right|\gs D\right) \\
&=O(1)\, \frac 1 d \sum_{j=1}^d\frac{\phi_i^2s_i^2}{(s_i^{2}+\phi_i^2)^2}=O(1)\,\frac 1 d A_d,
\end{align*}
which converges to 0 for a) and b) -- for c) we multiply the original statistic by $\sqrt{d}/A_d$, which means this term is multiplied with $d/A_d^2$ leaving us with $1/A_d$ which converges to 0 if $A_d/\sqrt{d}\to \infty$. Similarly, we can bound $\frac{1}{\sqrt{T}}\left( \sum_{t=1}^{\lfloor T x\rfloor}\eta_{d+1,t}(d)-\frac{\lfloor T x\rfloor}{T}\sum_{t=1}^T\eta_{d+1,t}(d) \right)$, showing that the middle term is asymptotically negligible. The assertions now follow by an application of the functional central limit theorem for\\ $\frac 1 T \left( \sum_{t=1}^{\lfloor T x\rfloor}\eta_{d+1,t}(d)-\frac{\lfloor T x\rfloor}{T}\sum_{t=1}^T\eta_{d+1,t}(d) \right)^2$.
\end{proof}

\begin{proof}[of Theorem~\ref{th_miss_cont}]
	The proof is analogous to the one of Theorem~\ref{th_ind_cont_2} on noting that $\HDE_3^2(\bD_d)=\frac{\sqrt{d}}{A_d}\,\HDE_2^2(\bD_d)$ and $\sigma_i^2=s_i^2+\Phi_i^2$ by using Lemma~\ref{th_panel_null_miss} c) above. Concerning the remainder term $\widetilde{R}_T(x)$ note that $e_{i,t}=s_i\eta_{i,t}+\Phi_i\eta_{d+1,t}$, so that the remainder term can be split into two terms. The first term can be dealt with analogously to the proof of Theorem~\ref{th_ind_cont_2} and is of order $O_P\left(\sqrt{\frac{1}{A_d}\,T\HDE_3(\bD_d)}\right)$, while for the second summand we get by an application of the Cauchy-Schwarz-inequality
	\begin{align*}
		&\sup_{0\ls x\ls 1}\left|\frac{1}{A_d}\sum_{i=1}^d\frac{\delta_i\phi_i}{\sigma_i^2}\left(\sum_{t=1}^{\lfloor Tx\rfloor}\eta_{d+1,t}-\frac{\lfloor Tx\rfloor}{T}\sum_{t=1}^T\eta_{d+1,t}\right)\right|
		=O_P(\sqrt{T})\,\sqrt{\frac{\sum_{i=1}^d\frac{\delta_i^2}{\sigma_i^2}}{A_d}}\\
		&=O\left( \sqrt{T\,\HDE_3^2(\bD_d)} \right).
	\end{align*}
\end{proof}

\begin{proof}[of Corollary~\ref{cor_miss_panel}]
	By an application of the Cauchy-Schwarz inequality it holds
\begin{align*}
	 &\bD_d^T\Lambda_d^{-1}\Sigma\Lambda_d^{-1}\bD_d=\sum_{i=1}^d\delta_{i,T}^2\frac{s_i^2}{(s_i^2+\Phi_i^2)^2}+\left(\sum_{i=1}^d\frac{\delta_{i,T}\Phi_i}{s_i^2+\Phi_i^2}\right)^2\\
	&\ls \sum_{i=1}^d\frac{\delta_{i,T}^2}{\sigma_i^2} \,\left( 1+\sum_{i=1}^d\frac{\Phi_i^2}{\sigma_i^2} \right)= \bD_d^T\Lambda_d^{-1}\bD_d \,(1+ A_d),
\end{align*}
which implies assertion a) on noting that
\begin{align*}
	\HDE_1^2(\bD_d,\qo)=\frac{(\bD_d^T\Lambda_d^{-1}\bD_d)^2}{\bD_d^T\Lambda_d^{-1}\Sigma\Lambda_d^{-1}\bD_d}.
\end{align*}
b)
		This follows immediately from Theorem~\ref{th_random_proj} since by $0<c\ls s_j^2\ls C<\infty$ as well as as $\Phi^2_i\ls C$, it follows that
		\begin{align*}
			\|\bD_d\|^2\sim \bD_d^T\;\mbox{diag}\left(\frac{1}{s_1^2+\Phi_1^2},\ldots,\frac{1}{s_d^2+\Phi_d^2}\right)\,\bD_d.
		\end{align*}
\end{proof}

\section*{Acknowledgements}
The first author was supported by the Engineering and Physical Sciences Research Council (UK) grants : EP/K021672/2 \& EP/N031938/1. Some of this work was done while the second author was at KIT where her position was financed by the Stifterverband f\"ur die Deutsche Wissenschaft by funds of the Claussen-Simon-trust. Furthermore, this work was supported by the Ministry of Science, Research and Arts, Baden-W\"{u}rttemberg, Germany. Finally, the authors would like to thank the Isaac Newton Institute for Mathematical Sciences, Cambridge, for support and hospitality during the programme 'Inference for Change-Point and Related Processes', where part of the work on this paper was undertaken.

\bibliographystyle{abbrvnat}
\bibliography{panel}

\end{document}